%% file: agt-3-9.tex
\documentclass{gtart}

\input agtout

\lognumber{9}
\volumenumber{3}
\volumeyear{2003}
\papernumber{9}
\published{26 February 2003}
\pagenumbers{207}{234}
\received{17 January 2003}
\accepted{3 February 2003}

\usepackage{amsmath, amssymb, graphicx, psfrag}

\def\psfraga <#1,#2> #3#4{\psfrag #3 
{\smash{\rlap{\kern #1 \raise #2\hbox{#4}}}}}

\def\RR{{\mathbb R}}

\def\HH{{\mathbb H}}

\def\ZZ{{\mathbb Z}}

\def\teich{\mathop{\rm Teich}}

\def\QF{{\cal Q}{\cal F}}

\def\ML{{\cal  ML}}

\def\PML{{\rm P}{\cal {ML}}}

\theoremstyle{plain}
\newtheorem{thm}{Theorem}[section]
\newtheorem{lemma}[thm]{Lemma}

\newtheorem{cor}[thm]{Corollary}
\newtheorem{prop}[thm]{Proposition}

\theoremstyle{definition}

\def\a{{\alpha}}
\def\b{{\beta}}
\def\d{{\delta}}

\def\g{{\gamma}}

\def\s{{\sigma}}

\def\th{{\theta}}

\def\dd{{\partial}}

\def\A{{\cal A}}
\def\B{{\cal B}}

\def\S{{\cal S}}
\def\L{{\cal L}}

\def\ch{\,{\rm cosh}\,}
\def\sh{\,{\rm sinh}\,}

\def\log{\,{\rm log}\,}

\begin{document}


\title[Limit points of lines of minima]{Limit points of lines of minima\\in Thurston's boundary of Teichm\"uller space}
\covertitle{Limit points of lines of minima\\in Thurston's boundary 
of Teichm\noexpand\"uller space}
\asciititle{Limit points of lines of minima\\in Thurston's boundary of 
Teichmueller space}

\author{Raquel D\'{\i}az\\Caroline Series}%
\coverauthors{Raquel D\noexpand\'{\noexpand\i}az\\Caroline Series}%
\asciiauthors{Raquel Diaz\\Caroline Series}%

\address{Deparmento Geometr\'{\i}a y Topolog\'{\i}a, 
Fac.\ CC.\ Matem\'aticas\\Universidad Complutense, 28040 Madrid, Spain}
\secondaddress{Mathematics Institute, University of Warwick\\Coventry 
CV4 7AL, UK}
\asciiaddress{Deparmento Geometria y Topologia, Fac. CC. 
Matematicas\\Universidad Complutense, 28040 Madrid, 
Spain\\and\\Mathematics 
Institute, University of Warwick\\Coventry 
CV4 7AL, UK}

\email{radiaz@eucmos.sim.ucm.es} 
\secondemail{cms@maths.warwick.ac.uk}
\asciiemail{radiaz@eucmos.sim.ucm.es, cms@maths.warwick.ac.uk}

\asciikeywords{Teichmueller space, Thurston  boundary, 
measured geodesic lamination,
Kerckhoff line of minima}

\begin{abstract}
Given two measured laminations $\mu$ and $\nu$ in a hyperbolic
surface which fill up the surface, Kerckhoff~\cite{KerckLM}  defines an
associated {\em line of minima} along which convex combinations of the
length functions of $\mu$ and $\nu$ are minimised. This is a line in
Teichm\"uller space which can
be thought as  analogous to the geodesic in hyperbolic space  determined
 by two points at
infinity. We show that when $\mu$ is uniquely
ergodic, this line converges to the projective lamination $[\mu]$, but
that when $\mu$ is rational,
the line converges not to $[\mu]$, but rather to the  barycentre of the
support of $\mu$.
Similar results on the behaviour of Teichm\"uller geodesics have been
proved by Masur~\cite{Masur}.
\end{abstract}

\asciiabstract{
Given two measured laminations mu and nu in a hyperbolic surface which
fill up the surface, Kerckhoff [Lines of Minima in Teichmueller space,
Duke Math J. 65 (1992) 187-213] defines an associated line of minima
along which convex combinations of the length functions of mu and nu
are minimised.  This is a line in Teichmueller space which can be
thought as analogous to the geodesic in hyperbolic space determined by
two points at infinity. We show that when mu is uniquely ergodic,
this line converges to the projective lamination [mu], but that
when mu is rational, the line converges not to [mu], but rather
to the barycentre of the support of mu.  Similar results on the
behaviour of Teichmueller geodesics have been proved by
Masur [Two boundaries of Teichmueller space, Duke Math. J. 49 (1982)
183-190].}

\maketitlepage


\section{Introduction}
\label{sec:introduction}

Let $S$ be a surface of hyperbolic type, and denote 
its Teichm\"uller
space by $\teich(S)$. Given a measured geodesic lamination $\mu$ on $S$
(see Section  \ref{sec:background} for definitions), there is a function
$l_{\mu} \co \teich(S) \to \RR ^+$ which associates to each  $\rho\in
\teich(S)$  the hyperbolic length $l_{\mu}(\rho)$ of $\mu$ in the
hyperbolic
structure $\rho$.  If $\mu, \nu$ are two measured laminations which fill
up the surface, Kerckhoff \cite{KerckLM} proved that for any number
$s\in (0,1)$, the function $F_s= (1-s) l_{\mu} +s l_{\nu}$ has a 
global minimum at a unique point $m_s \in \teich(S)$. The set of all
these minima, when $s$ varies in the
interval $(0,1)$, is called a \emph{line of minima} $\L_{\mu,\nu}$.

The Teichm\"uller space of a surface is topologically a ball which, as
shown by Thurston, can be
compactified by the space $\PML$ of projective measured laminations on
$S$. Various
analogies between Teichm\"uller space and hyperbolic space have been
studied, for example {\it earthquake paths} in Teichm\"uller space are
analogous to horocycles in hyperbolic space.
In  \cite{KerckLM}, Kerckhoff studied some properties of the
lines of  minima  analogous to properties of geodesics in hyperbolic
space. For example, two projective measured laminations determine a
unique line of minima, in analogy to the fact that two different
 points in the boundary of hyperbolic space determine a unique geodesic.
He warns, however,  that lines of minima
do not always converge  to the point corresponding to $\mu$
in Thurston's compactification of $\teich (S)$, mentioning that examples
can be constructed by taking $\mu$ rational (that is, such that its
support consists entirely of closed leaves).  In this paper we make this
explicit by showing that any line of minima $\L_{\mu,\nu}$ where
$\mu=\sum a_i\a_i$, with $a_i > 0$,  converges to the
projective lamination $[ \sum \a_i ]$, rather than to $[ \sum a_i \a_i ]$.

\begin{thm}
\label{thm:linetobarycenter}
Let $\mu=\sum_{i=1}^N a_i \a_i$  be a rational measured lamination
(that is,  $\a_i$ is a collection
of disjoint simple closed curves on $S$ and $a_i> 0$ for all $i$)
 and $\nu$ any measured lamination so that $\mu, \nu$   fill up
 the surface. 
For any  $0<s<1$, consider the function $F_s \co \teich(S) \to \RR $
defined by $F_s(\rho)=(1-s) l_{\mu}(\rho)+s l_{\nu}(\rho)$, and denote
 its unique minimum  by
$m_s$. Then
\[ \lim_{s\to 0} m_s=[\a_1+ \dots + \a_N]  \in \PML. \]
\end{thm}

By contrast, if $\mu$ is uniquely ergodic and maximal (see
Section~\ref{sec:background} for the
definition), we prove:

\begin{thm}
\label{thm:uniqueergodic}
Let $\mu$  and $\nu$ be two measured laminations which fill up the
surface and such that $\mu$ is uniquely ergodic and maximal.
With $m_s$  as above,  
$$ \lim_{s\to 0} m_s = [\mu]  \in \PML.$$
\end{thm}

Exactly similar results have been proved by Masur~\cite{Masur} for
Teichm\"uller geodesics. In this case, a geodesic  ray is determined by
a base surface $\rho $ and a quadratic differential $\phi$ on $\rho$.
Roughly
speaking, the end of this ray depends on the horizontal foliation $F$ of
$\phi$. Masur  shows that if $F$ is a Jenkins--Strebel differential, that
is, if 
 its horizontal foliation has closed leaves,  then the associated
ray converges in the Thurston boundary to the barycentre of the leaves
(the foliation with the same closed leaves all of whose cylinders have
unit height), while if $F$ is  uniquely ergodic and every leaf (apart
from saddle connections) is dense in $S$, it converges to the
boundary point defined by $F$. 
 
Our interest in lines of minima arose from the study of the space $\QF
(S)$ of quasifuchsian groups associated to a surface $S$. The {\it
pleating plane} determined by a pair of projective measured
laminations is the set of quasifuchsian groups whose convex hull
boundary is bent along the given laminations with bending measure in
the given classes.  It is shown in~\cite{SKerck}, see also
\cite{SKerck2}, that if $\mu, \nu$ are measured laminations, then the
closure of their pleating plane meets fuchsian space exactly in the
line of minima $\L_{\mu,\nu}$.
 
From this point of view,  it is often more natural to look at the
collection of all groups whose  convex hulls are bent along a specified
set of closed curves. That is, we  forget  the proportions between the
bending angles given by the measured lamination $\mu$ and look  only at
its support. 
This led us in~\cite{DS}  to study the 
{\em simplex of minima}   determined by two 
systems of  disjoint simple
curves on the twice punctured
torus, where  
direct calculation of some special examples led to our results here.

The {\it simplex of minima} $\S_{\A,\B}$ associated to  systems of
disjoint simple curves 
$\A=\{\alpha_1, \ldots, \alpha_N\}$ and 
$\B=\{\beta_1,\ldots,\beta_M\}$ which  fill  up  the
surface,  is the union of lines
of minima ${\cal L}_{\mu,\nu}$, where $\mu,  \nu\in \ML (S)$ are
strictly positive 
linear combinations of $\{\alpha_i\}$ and  $\{\beta_i\}$, respectively.  
We can regard  $\S_{\A,\B}$ as the image of the  affine
simplex $S_{\A,\B}$ in  $\RR^{N+M-1}$ spanned by independent points  
$A_1,\dots,A_N,B_1,\dots, B_M$, under the map $\Phi$  which sends the
 point
$(1-s)(\sum_i a_iA_i)+s(\sum_jb_jB_j)$  (with $0< s,a_i,b_j< 1, \sum
a_i=1,\sum b_j=1$) to  the unique minimum of the function $(1-s)(\sum_i
a_il_{\a_i})+s(\sum_j b_jl_{\b_j})$.

As observed in~\cite{DS}, the methods of~\cite{KerckLM} show that the
map $\Phi$ is continuous and proper. It may or may not be a
homeomorphism onto its image;  in \cite{DS} we give a necessary and
sufficient condition and show by example that both cases occur. The map
extends continuously to the faces of $S_{\A,\B}$
which correspond  to curves $\{\a_{i_1},\ldots,\a_{i_k}\}$,
$\{\b_{j_1},\dots,\b_{j_l}\}$ that still fill up the surface.
Nevertheless, as a consequence of Theorem~\ref{thm:linetobarycenter}, 
$\Phi$ does not necessarily extend to a function from the closure of  
$S_{\A,\B}$ into the Thurston boundary. 

\begin{cor}
\label{cor:nonextend} Let $\A, \B$ be as above and suppose that
$\{\alpha_1, \ldots, \alpha_{N-1}\}$ and 
$\B=\{\beta_1,\ldots,\beta_M\}$ also fill up $S$.
Then, the map $\Phi\co  S_{\A,\B} \to \teich(S)$ does not extend
continuously to a function  $\overline{ S_{\A,\B}}\to \teich(S)\cup
\PML(S)$.
\end{cor}
\begin{proof}
 Let  $\{x_n\}$  be a sequence of points in $S_{\A,\B}$, and $  \{y_n\}$
another sequence  in the face spanned by $A_1,\dots,A_{N-1}, B_1,\dots,
B_M$, both converging to  $(A_1+\dots +A_{N-1})/(N-1)$. Then, by
Theorem~\ref{thm:linetobarycenter},  $\Phi(x_n)$ converges to
$[\a_1+\dots+\a_{N}]$ while $\Phi(y_n)$ converges to
$[\a_1+\dots+\a_{N-1}]$.
\end{proof}

We remark that examples of curve systems as in the corollary are easy to
construct. 

\medskip

The paper is organised as follows. The main work is in proving  
Theorem~\ref{thm:linetobarycenter}.  In Section 2
 we recall background and give the (easy) proof of
Theorem~\ref{thm:uniqueergodic}.  In  Section 3 we study
an example which illustrates Theorem~\ref{thm:linetobarycenter} and its
proof. The general proof is easier when $\a_1,\dots , \a_N$ is
a pants decomposition. We work this case in Sections 4 and 5,
leaving the non-pants decomposition case for Section 6.

\medskip
The first author would like to acknowledge  
partial  support from  MCYT grant BFM2000-0621 and 
UCM grant PR52/00-8862, and the second
support from an EPSRC Senior Research Fellowship.

\section{Background}
\label{sec:background}

We take the Teichm\"uller space $\teich(S)$ of a surface $S$ of
hyperbolic type to be the set
of faithful and discrete representations 
$\rho \co \pi_1(S) \to PSL(2,\RR)$ which take loops around punctures
to
parabolic elements, up to conjugation by elements of $PSL(2,\RR)$. An
element $\rho$ of $\teich(S)$ can be regarded as a {\it marked}
hyperbolic structure on $S$. The space $\teich(S)$ is topologically a
ball of dimension $2(3g-3+b) $, where $g$ is the genus and $b$ the
number of punctures of $S$. A {\it pants decomposition} of $S$ is a set
of disjoint simple closed curves, $\{\a_1,\dots,\a_N\}$ which decompose
the surface into pairs of pants ($N=3g-3+b$).  Given  a pants
decomposition $\{\a_1,\dots,\a_N\}$ on $S$, the {\it Fenchel-Nielsen}
coordinates give a global parameterization of $\teich(S)$.  Given a
marked
hyperbolic structure on $S$, these coordinates consist of the lengths  
$l_{\a_i}$ of the geodesics representing  the curves $\a_i$, and the
{\it twist
parameters} $t_{\a_i}$. The lengths $l_{\a_i} $   determine uniquely the
geometry on each pair of pants, while the twist parameters 
 are real numbers  determining the way these pairs of pants are glued
together to build up the hyperbolic surface. We need to specify a set of 
base points in $\teich(S)$, namely a subset of $\teich(S)$  where the
twist parameters are all equal to zero. This can be done by choosing a
set of curves $\{\d_i\}$ {\it dual } to the $\{\a_i\}$, in the sense
that each $\d_i$ intersects $\a_i$ either  once or twice and is disjoint
from $\a_j$ for all $j\not=i$. For each fixed set of values of
$l_{\a_i}$, the
base point is then the marked hyperbolic structure in which each $\d_i$
is orthogonal to $\a_i$, when they intersect once, or in which the two
intersection angles (measured from $\d_i$ to $\a_i$)   sum to $\pi$,
when they intersect twice. 

A  {\it geodesic lamination} in a hyperbolic surface $\rho$ is a closed
subset of $\rho$ which is disjoint union of   simple geodesics, called
its {\it leaves}. 
A geodesic lamination is {\it measured} if it carries a transverse
invariant measure (for details, see for example~\cite{FLP, Penner} and
the appendix to~\cite{Otal}). 
The space $\ML$ of measured laminations is
given the weak topology. If $\mu \in \ML$, then $|\mu|$ will denote its
underlying support. To exclude trivial cases, we assume that each leaf
$l$ of $|\mu|$ is a density point of $\mu$, meaning that any open
interval transverse to $l$ has positive $\mu$-measure. 
 In this paper, we shall mainly use {\it rational } measured
laminations, denoted by $\sum_i a_i\a_i $, where $a_i \in \RR^{+}$ and 
$\a_i$ are disjoint simple closed geodesics. This measured lamination
assigns mass $a_i$ to each intersection of a transverse arc with $\a_i$.
The length of a rational lamination 
$\sum_i a_i\a_i $ on a hyperbolic surface $\rho$  is defined to be
$\sum_i a_i l_{\a_i}(\rho)$, where $l_{\a_i}(\rho)$ is the length of
$\a_i$ at $\rho$.   Rational measured laminations are dense in $\ML$ and  
the length of a measured lamination can be defined as the limit of the
lengths of approximating rational measured laminations,
see~\cite{KerckEA}.  This construction appears to depend on $\rho$,
however a homeomorphism  
between   hyperbolic surfaces  transfers geodesic
laminations  canonically  from the first surface to the second. Thus,
given a
measured lamination $\mu$,  there is a map $l_{\mu} \co \teich(S)
\to \RR^+$ which assigns to each point $\rho\in \teich(S)$ the length
$l_{\mu}(\rho)$ of $\mu$ on the hyperbolic structure $\rho$.  The map
$l_{\mu}$ is real analytic with respect to the real analytic structure
of $\teich(S)$, see~\cite{KerckEA} Corollary 2.2. 

Two measured laminations are {\em equivalent} if they have the same
underlying
support and proportional transverse measures. The equivalence class of a
measured lamination $\mu$ is called a {\it projective measured
lamination } and is denoted by $[\mu]$. A measured lamination is {\it
maximal}
 if its support is not contained in the support of any other measured
lamination. A lamination is {\em uniquely ergodic} if every lamination
with the same support is in the same projective equivalence class. (Thus
the lamination 
$\sum_i a_i {\a_i}$ is uniquely ergodic if and only if the sum contains
exactly one term.)
The geometric intersection number 
$i(\g,\g')$ of two simple closed geodesics is the number of points in
their intersection. 
This number  extends by  bilinearity and continuity to the intersection
number of measured laminations, see~\cite{ Rees, KerckEA, Bon}. 
The following characterisation of uniquely ergodic laminations is needed
in the proof of Theorem 1.2.

 \begin{lemma}
\label{lemma:charue} A lamination $\mu \in \ML$ is uniquely ergodic and
maximal if and only if,  for  all $\nu \in \ML$,   $i(\mu,\nu) =0$
implies $\nu \in [\mu]$.
\end{lemma}
\begin{proof} If $i(\mu,\nu) =0$ implies $\nu \in [\mu]$, it is easy to
see that $\mu$ must be uniquely ergodic and maximal.
The converse  follows using the definition of intersection
number as the integral over $S$ of the product measure
$\mu \times \nu$; see for example~\cite{KerckEA}.
Since we are assuming $\mu$ is uniquely ergodic, it is enough to show
that the supports of $\mu$ and $\nu$ are the same.
 
Let $\omega$ be the lamination consisting of leaves (if any) which are
common to $|\mu| $ and $ |\nu|$. Let $\mu_{\omega}$ and $
\nu_{\omega}$ denote the restrictions of $\mu$ and $\nu$ to $\omega$.
Clearly $\omega$ is closed, and hence (using the decomposition of
laminations into finitely many minimal components, see for
example~\cite{CasBl},\cite{Otal}), one can write $\mu = \mu_{\omega} +
\mu'$, $\nu = \nu_{\omega} + \nu'$ where $\mu',\nu'$ are (measured)
laminations disjoint from $\omega$ such that every leaf of $\mu'$ is
transverse to every leaf of $\nu'$. Since $\mu$ is uniquely ergodic,
one or other of $\mu'$ or $\mu_{\omega}$ is zero. In the former case,
maximality of $\mu$ forces $\nu'=0$, and we are done.

Thus we may assume that every leaf of $|\mu |$ intersects every leaf of 
$ |\nu|$ transversally; let $X$ be the   set of intersection points of 
these leaves. Cover $X$ by small disjoint open `rectangles'
$R_i$,  each with two `horizontal' and two `vertical' sides, in such a
way that $|\mu| \cap R$ consists entirely of arcs with endpoints on the
horizontal sides and similarly for  $|\nu| \cap R$ replacing horizontal
by vertical.
Put a product measure on $R$ by using the transverse measure $\mu$ on
`horizontal' arcs and $\nu$ on `vertical' ones. Then  $i(\mu,\nu) =
\sum_i \int_{R_{i}} d\mu \times d \nu$.  Our assumption that each leaf
of $|\mu|$ and $|\nu|$ is a density point means that the contribution to
$i(\mu,\nu)$ is non-zero whenever $R \cap X$ is non-empty. 
Thus $i(\mu,\nu)=0$ implies that $X=\emptyset$. Since $\mu$ is maximal,
every leaf of $|\nu|$ either coincides with or intersects some leaf of
$|\mu|$, and we conclude that the  leaves of $|\nu|$ and  $|\mu|$
coincide as before.
\end{proof}

There is a similar characterisation of uniquely ergodic foliations due
to Rees~\cite{Rees}, see also~\cite{Masur} Lemma 2, in which the
assumption that $\mu$ is maximal is replaced by the assumption that
every leaf, other than saddle connections, is dense. 
(Notice  that the above proof shows that a uniquely ergodic lamination
is also minimal, in the sense that every leaf is dense in the whole
lamination.)

\subsection{The Thurston Boundary} 
 We denote the set of all  non-zero projective measured laminations  
on $S$ by 
 $\PML(S)$.
Thurston has shown that
$\PML(S)$ compactifies $\teich(S)$ so that $\teich(S)\cup \PML(S)$ is
homeomorphic to a closed ball. We explain this briefly; for details
see~\cite{FLP}. A sequence $\{ \rho_n\}\subset \teich(S)$ converges to 
$[\mu] \in \PML$ if the lengths of simple closed
geodesics on $\rho_n$ converge projectively to their intersection
numbers with $\mu$;  more precisely, if there exists a sequence
$\{c_n\}$ converging to infinity, so that 
$ l_{\g}(\rho_n)/c_n \to i(\g,\mu)$, for any simple closed geodesic
$\g$.
The following lemma summarises the consequences of  this definition we
shall need.

\begin{lemma}
\label{lemma:limit}
Let $\a_1,\dots,\a_N$ be a pants decomposition on $S$ and 
let $\{\rho_n\} \subset \teich(S)$ so that $\rho_n\to [\mu]\in \PML(S)$.
Then:\begin{itemize}
\item[\rm(a)] if $\nu\in \ML$ with  $i(\mu,\nu)\not=0$ then
$l_{\nu}(\rho_n)\to \infty$,
\item[\rm(b)] if  $l_{\a_i}(\rho_n)$ is bounded for all $i=1,\dots,N$,
then there exist $a_1,\dots,a_N \geq 0$ so that $[\mu]= [a_1\a_1+\dots+
a_N\a_N]$.\end{itemize}
\end{lemma}

The proofs are immediate from the definitions. Part  (b) gives a 
 sufficient condition for convergence to a rational lamination
$[\Sigma_i a_i\a_i]$, when the $\{\a_i\}$ is a pants decomposition. To
compute the coefficients $a_i$ we  take another  system of curves
$\{\d_i\}$ dual to the $\{a_i\}$. From the definition, 
$$ {{l_{\d_j}(\rho_n)}\over{l_{\d_k}(\rho_n)}}\to{{i(\Sigma
a_i\a_i,\d_j) }\over{i(\Sigma
a_i\a_i,\d_k})}={{a_ji(\a_j,\d_j)}\over{a_ki(\a_k,\d_k)}},
$$
and we know that $i(\a_i,\d_i)$ is either 1 or 2, so this gives the
proportion $a_j/a_k$.

\subsection{Lines of minima} 
Two measured laminations $\mu,\nu$ are said to {\em fill up} a surface
$S$ if for any other
lamination $\eta$ we have $i(\mu,\eta)+i(\nu, \eta)\not=0$. It is proved
in~\cite{KerckLM}  that for any two such
laminations,  the function  $ l_{\mu}(\rho)+ l_{\nu}(\rho)$
has a unique minimum on $\teich(S)$.
Thus $\mu$ and $\nu$
determine the {\em line of minima} $\L_{\mu,\nu}$, namely the
set of points $m_s\in \teich(S)$ 
at which the function  $F_s(\rho)=(1-s)l_{\mu}(\rho)+s l_{\nu}(\rho)$
reaches its minimum as $s$ varies in $  (0,1)$.

Given this definition, we can immediately prove
Theorem~\ref{thm:uniqueergodic}.

\begin{proof}[Proof of Theorem~\ref{thm:uniqueergodic}]
Observe that $l_{\mu}(m_s)$ is bounded as $ s \to 0$,
because $ l_{\mu}  \le 2((1-s)   l_{\mu} +s   l_{\nu}) = 2F_s $ for $s <
1/2$ 
and  $F_s(m_s) \le F_s(\rho_0) $ where $\rho_0$ is some arbitrary  point
in $\teich(S)$. By compactness of $\teich (S) \cup \PML$, we can  choose
some sequence $s_n \to 0$ for which $m_{s_n}$ is convergent. Moreover,
it is proved in \cite{KerckLM} that  the map $s \to m_s$ is proper, and
so   $m_{s_n} \to [\eta] \in \PML$. By Lemma~\ref{lemma:limit} (a) we
have that 
$i(\mu,\eta) = 0$ and from  Lemma~\ref{lemma:charue} we deduce that $
[\eta] = [\mu]$. The result follows.
\end{proof}

We now turn to  the more interesting rational case.
In general, the minimum $m_s$ is in fact the unique critical point of
$F_s$, so 
a point $p \in \L_{\mu,\nu}$ if and only if 
the  1-form ${\rm d} F_s =(1-s) {\rm d} l_{\mu} +s {\rm d} l_{\nu}$
vanishes at $p$ for some $s$.  If
$\mu=\sum_i^N a_i\a_i$ is rational and  $\{\a_1,\dots ,
\a_N\}$ is a pants decomposition, this enables us to
find equations for  $\L_{\mu,\nu}$. In fact, applying 
 ${\rm d} F_s = (1-s)\sum_i a_i {\rm d} l_{\a_i} +s {\rm d} l_{\nu}$  to
the tangent vectors $\frac{\partial}{\partial t_{\a_i}}$, we get 
 \begin{equation}
\label{eq:form1}
\frac{\partial l_{\nu}}{\partial t_{\a_i}}=0, \quad  \hbox{for all}
\quad i=1,\dots,n.
\end{equation}
 Similarly, applying ${\rm d} F_s$ to the tangent vectors 
$\frac{\partial}{\partial l_{\a_i}}$, we get 
 $
{\partial l_{\nu}}/{\partial l_{\a_i}}=-a_i (1-s)/{s},
$ 
 so that the line of minima satisfies the equations
\begin{equation}
\label{eq:form2}
\frac{1}{a_i}\frac{\partial l_{\nu}}{\partial
l_{\a_i}}=\frac{1}{a_j}\frac{\partial l_{\nu}}{\partial l_{\a_j}}\quad 
\hbox{for all} \quad i,j.
\end{equation} 
Since  $\frac{\partial}{\partial l_{\a_i}}$, $\frac{\partial}{\partial
t_{\a_i}}$ form a basis of tangent vectors (\cite{KerckEA} Proposition
2.6), the equations 
(\ref{eq:form1}) and (\ref{eq:form2}) completely determine
$\L_{\mu,\nu}$.

\section{Example}
\label{sec:example}

Let $S=S_{1,2}$ be the twice punctured torus. 
Consider   two
disjoint, non-disconnecting simple closed curves $\a_1$ and $\a_2$, and
let $\b$ be a simple
closed curve intersecting each of $\a_1$ and $\a_2$ once. 
 For positive numbers $a_1,a_2$, denote by $\mu$ the measured
lamination $a_1\a_1+a_2\a_2$. 
(When $S$ has a
hyperbolic structure and $\gamma  \in \pi_1(S)$,  we abuse  notation by
using
$\gamma$ to mean also the unique
geodesic  in the homotopy class of $\gamma$.)
We shall compute the equation of  the line
of minima $\L_{\mu,\b}$, in terms of Fenchel-Nielsen coordinates
relative to the pants decomposition $\{\a_1,\a_2\}$ and  dual curves
$\{\d_1,\d_2\}$ (see Figure 1),
  and we shall
show that this line converges in Thurston's compactification to
$[\a_1+\a_2]$.

As explained above, the line of minima  $\L_{\mu,\b}$ is determined by
the equations 
$$ \frac{\dd l_{\b}}{\dd t_{\a_1}}=0,  \quad \frac{\dd l_{\b}}{\dd
t_{\a_2}}=0
  \quad {\rm and} \quad 
 \frac{1}{a_1}\frac{\dd l_{\b}}{\dd l_{\a_1}}=\frac{1}{a_2} \frac{\dd
l_{\b}}{\dd l_{\a_2}}.
$$
\begin{figure}[ht!]\small
\psfrag{a_1}{$\alpha_1$}
\psfrag{a_2}{$\alpha_2$} 
\psfrag{b}{$\beta$}
\psfrag{d_1}{$\delta_1$}
\psfrag{d_2}{$\delta_2$}
\centerline{\includegraphics[height=5cm]{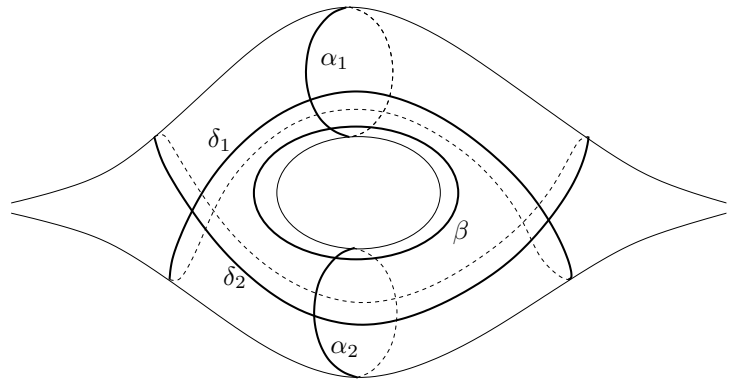}}
\caption{The curves $\a_1,\a_2$ and $\b$ on a twice punctured torus}
\end{figure}%
For two simple closed geodesics $\a,\b$,  Kerckhoff's derivative
formula~\cite{KerckNRP} states that $(\partial l_{\b}/\partial
t_{\a})|_{\rho}= \sum \cos
\psi_i(\rho)$, where $\psi_i$ are the intersection angles from $\b$ to
$\a$ at each intersection  point. Thus the first two equations mean 
that 
at a point in the line of minima the geodesic $\b$ is orthogonal to
$\a_1$
and $\a_2$. Let $P, P'$ be the two pairs of pants into which
$\{\a_1,\a_2\}$ split $S$; denote by $H_{12},H_{11}$ the perpendicular
segments  in $P$ from the geodesic $\a_1$ to $\a_2$ and from $\a_1$ to
itself, respectively; and denote by  $H_{12}',H_{11}'$ the analogous
perpendiculars in $P'$. Since $\b$ is orthogonal to $\a_1$ and $\a_2$,
$P$ and $P'$ are glued so that the segments $H_{12}$ and $H_{12}'$ match
up. 
 Since $P, P'$  are isometric (each is  determined by the
lengths $(l_{\a_1},l_{\a_2},0)$), the segments $H_{11}$ and $H_{11}'$
also match, so that the union of both segments is the geodesic 
$\d_1$. Therefore, $\d_1$ intersects $\a_1$ orthogonally, and so the
twist parameter $t_{\a_1}$ is zero. In the same way, $\d_2$ intersects
$\a_2$ orthogonally and $t_{\a_2}=0$.

It is not difficult to find the expression for the
length of $\b$ in the Fenchel-Nielsen coordinates
$(l_{\a_1},l_{\a_2},t_{\a_1},t_{\a_2})$, either using trigonometry  or
by looking at the trace of the
element representing $\b$ in the corresponding fuchsian group.
This is done in detail in~\cite{DS}. We have
$$
\ch \frac{l_{\b}}{2}= \frac{1+ \ch \frac{l_{\a_1}}{2}
\ch\frac{l_{\a_2}}{2}} {\sh\frac{l_{\a_1}}{2}\sh\frac{l_{\a_2}}{2}}   
\ch \frac{t_{\a_1}}{2} \ch \frac{t_{\a_2}}{2}+
\sh  \frac{t_{\a_1}}{2} \sh \frac{t_{\a_2}}{2}.
$$ 
Computing the  derivatives $\partial l_{\b}/ \partial l_{\a_i}$
directly from this formula   we get
$$
-\sh \frac{l_{\b}}{2} \frac{\dd l_{\b}}{\dd l_{\a_1}} =\frac{\ch
\frac{l_{\a_1}}{2} + \ch\frac {l_{\a_2}}{2} } {\sh^2 \frac{l_{\a_1}}{2} 
\sh\frac{l_{\a_2}}{2} }    \ch \frac{t_{\a_1}}{2} \ch\frac{t_{\a_2}}{2}
,
$$
$$ 
-\sh \frac{l_{\b}}{2} \frac{\dd l_{\b}}{\dd l_{\a_2}} =\frac{
\ch \frac{l_{\a_1}}{2} + \ch \frac{l_{\a_2}}{2} } {\sh
\frac{l_{\a_1}}{2}  \sh^2 \frac{l_{\a_2}}{2} }   \ch \frac{t_{\a_1}}{2}
\ch\frac{t_{\a_2}}{2}  .
$$ 
Therefore,  the  equations determining the line of minima are: 
$$  t_{\a_1}=t_{\a_2}=0, \quad \frac{a_1}{a_2} = \frac{\sh (l_{\a_2}/2)
}{\sh(l_{\a_1}/2) }. 
$$
By allowing $a_1, a_2 $ to vary among all positive
numbers, we observe that the corresponding lines of minima are pairwise
disjoint and in fact foliate  the whole plane $ \{
(l_{\a_1},l_{\a_2},t_{\a_1},t_{\a_2}) \;| \;
t_{\a_1}=0,t_{\a_2}=0 \}$
in $\teich(S_{1,2})$.

Clearly, at one end of $\L_{\mu,\b}$  the lengths $l_{\a_1},
l_{\a_2}$ tend to zero. This cannot happen when $s\to 1$ because if
$l_{\a_1},
l_{\a_2}\to 0$, then   $l_{\b}$  and hence 
$(1-s)l_{\mu}(m_s)+sl_{\b}(m_s)$ tend to $\infty$, and  this contradicts 
the fact that $m_s$ is the minimum.  
  Thus  $l_{\a_1}, l_{\a_2}\to 0$ as  $s \to 0$ and therefore by
Lemma~\ref{lemma:limit} (b),  the line $\L_{\mu,\b}$
converges to a point of the form $[a'_1\a_1+a'_2\a_2]$, for some
$a'_1,a'_2 \geq 0$. To compute these numbers, we compute the lengths of
the dual  curves $\d_1,\d_2$. By hyperbolic trigonometry  we get
$$ \ch\frac{ l_{\d_1}}{4} =  \sh \frac{l_{\a_2}}{2}\sh \frac{l_{\b}}{2},
\quad
\ch\frac{ l_{\d_2}}{4} =  \sh \frac{l_{\a_1}}{2}\sh \frac{l_{\b}}{2}.
 $$
Thus comparing the lengths of $\d_1,\d_2$ along the line of minima
$\L_{\mu,\b}$ we find
$$
\lim_{s\to 0} \frac{\ch( l_{\d_1}/4)}{\ch (l_{\d_2}/4)} =\lim _{s \to 0}
\frac{\sh(l_{\a_2}/2)}{\sh(l_{\a_1}/2)} =\frac{a_1}{a_2}.
$$
Since $l_{\a_i}\to 0$, we have that $l_{\d_i}\to \infty$, so that
 $$ 
\frac{a_1}{a_2}=
\lim_{s\to 0} e^{(l_{\d_1}-l_{\d_2})/2} .$$
Taking logarithms, we get that $\lim_{s\to 0}(l_{\d_1}-l_{\d_2})$ is a
constant, and this implies that  $\lim_{s\to 0}(l_{\d_1}/l_{\d_2})=1$.
Hence,  $\L_{\mu,\b} \to [\a_1+\a_2]$ as $s\to 0$.

\section{Statements of main results: pants decomposition case}
\label{sec:statements}

In order to prove Theorem~\ref{thm:linetobarycenter}, we need to show
that the lengths of all the simple closed geodesics converge
projectively to their intersection numbers with $\sum \a_i$.  So we
want to estimate the length of any simple closed geodesic along the
line of minima.  We first prove that along the line of minima the
lengths of $\a_i$ tend to zero and the twist parameters about $\a_i$
are bounded (Proposition~\ref{prop:geom.prop.LM} (a) and (b)). When
$\a_i$ is a pants decomposition, these two properties allow one to
give a nice estimate of the length of a closed geodesic
(Proposition~\ref{prop:lengthestimate}): the main contribution is
given by the arcs going through the collars around the curves
$\a_i$. Finally, to compare the length of two closed geodesics, we
need to compare the orders of the lengths $l_{\a_i},l_{\a_j}$. This is
done in Proposition~\ref{prop:geom.prop.LM} (c).

  In this section we state these two propositions in the case in which
$|\mu|$ is a pants decomposition and $\nu$ is rational.    The
propositions  will  be proved in the next section.
  Both propositions remain true when   $\nu$ (in
Proposition~\ref{prop:geom.prop.LM})  and $\g$ (in
Proposition~\ref{prop:lengthestimate}) are arbitrary  measured
laminations.
We will comment on the proof of these stronger versions in
Section~\ref{sec:irrational}.

 Recall that two real functions $f(s), g(s)$ have the {\it same
order} as $s\to s_0$, denoted by $f\sim g$, when there exist positive
constants $k<K$ so that $k< f(s)/g(s)<K$ for all $s$ near enough  to
$s_0$. Write $f(s) = O(1)$ if $f(s)$ is bounded. 
\label{page:order}

\begin{prop}
 \label{prop:geom.prop.LM}
Suppose that  $\mu=\sum a_i\a_i$ and $\nu= \sum b_i\b_i$ are two
measured
laminations where $\{\a_i\}$ is a pants decomposition, and $a_i>0$ for
all $i$. Let $m_s$ be the minimum point of the function $F_s$. Then 
\begin{itemize}
\item[\rm(a)] for any $i$, $\lim_{s\to 0}l_{\a_i}(m_s)=0$;
\item[\rm(b)] for any $i$, $|t_{\a_i}(m_s)|$ is bounded when $s\to 0$;
\item[\rm(c)] for all $i,j$, $l_{\a_i}(m_s)\sim l_{\a_j}(m_s)\sim s$ as
$s\to 0$.\end{itemize}
\end{prop}
The proof of Proposition~\ref{prop:geom.prop.LM} (a) is direct and could
be read now.

\begin{prop}  
\label{prop:lengthestimate} Let $\{\a_1,\dots,\a_N\}$ be a pants
decomposition of $S$ and $\g$ a closed geodesic.
Let $\rho_n\in \teich(S)$ be  a sequence so that, when $n \to \infty$, all 
the lengths 
  $l_{\a_i}(\rho_n)$ are bounded above and the twists 
$t_{\a_i}(\rho_n)$  are bounded for all $i$.  Then, as $n \to \infty$,
we have
$$l_{\g}(\rho_n)=2\sum_{j=1}^N i(\a_j,\g) \log
\frac{1}{l_{\a_j}(\rho_n)} +O(1).$$
\end{prop}

 In view of this proposition, it is enough to work with collars around
$\a_j$ of width $2\log(1/l_{\a_j})$, even if they are not the maximal
embedded collars.
 The more relaxed hypothesis about the lengths $l_{\a_i}$
being bounded above is not needed if $\{\a_1,\dots,\a_N\}$ is a  pants 
decomposition, but
will be useful in the general case in Section 6.

\begin{proof}[Proof of Theorem~\ref{thm:linetobarycenter} for the pants
decomposition case] \nl
Suppose $\a_1,\dots,\a_N$ is a pants decomposition
system. 
By Proposition~\ref{prop:geom.prop.LM} (a), the lengths $l_{\a_i}$ tend
to zero as $s\to 0$. Therefore, by Lemma~\ref{lemma:limit} (b), $ m_s
\to
[a_1'\a_1+\dots +a_N'\a_N]$, as $s\to 0$, for some $a_i'\geq 0$. By
Proposition~\ref{prop:geom.prop.LM} (b), along the line of minima the
twists $t_{\a_i}$ are bounded. Then we can use
Proposition~\ref{prop:lengthestimate} to estimate the length of two
simple closed curves $\g,\g'$: the  proportion between their lengths is
\begin{equation}
\label{eq:lengthestimate}
\frac{l_{\g}}{l_{\g'}}=\frac{2\sum i(\a_j,\g) \log
\frac{1}{l_{\a_j}}+O(1)}{2\sum i(\a_j,\g') \log
\frac{1}{l_{\a_j}}+O(1)}. 
\end{equation}
Now, by Proposition~\ref{prop:geom.prop.LM}(c),   $l_{\a_i}\sim
l_{\a_j}$ as $s\to 0$; this implies that 
$  \log \frac{1}{l_{\a_i}}/\log \frac{1}{l_{\a_j}}$ $\to 1$ as $s\to 0$
(see Lemma~\ref{lemma:order} below).
 Dividing numerator and denominator of (\ref{eq:lengthestimate}) by
$\log (1/l_{\a_1})$, we get  that
$\lim_{s\to 0} (l_{\g}/l_{\g'}) =i(\sum \a_j,\g)/i(\sum \a_j, \g')$.
Hence $a'_i=1$ for all $i$.
\end{proof}

\section{Proof: pants decomposition case}
\label{sec:proofs}

 We shall estimate the length of a geodesic $\g$ comparing it with the
length of a ``broken arc" relative to a pants decomposition. Broken arcs
are a main tool in \cite{SW}, and we refer there for details. 
The idea is that there is a unique curve freely homotopic to $\g$ with
no
backtracking,  made up of arcs which wrap around a pants curve,
alternating with  arcs which cross a pair of pants from one boundary to
another following the common perpendiculars between the boundaries. This
collection of mutually perpendicular arcs constitute the broken arc,
whose length, as shown in
Lemma~\ref{lemma:brokenarc}, approximates the length of $\g$.   

 To determine the length of the broken arc, we study  the
geometry of a pair of pants. By using the trigonometric formulae for
right angle hexagons and pentagons, we can compute the length of the
segments perpendicular to two boundary components, and estimate this
length when the lengths of the boundary components tend to zero. This is
done in Lemma~\ref{lemma:pairofpants}.

\subsection{Broken arcs }
\label{sec:brokenarcs}

A {\it  broken arc}    
 in $\HH^2$ is a sequence of oriented segments such that
the final point of one segment is equal to the initial point of the
next,
 and such that consecutive arcs meet orthogonally.
Labelling the segments in order $V_1,H_1,\ldots, V_r,H_r,V_{r+1}$,  we
also
 require that for   $ 1 \le i \le 
r-1$  the segments $H_i,H_{i+1}$ are contained in opposite
halfplanes with respect to  $V_{i+1}$. We call  the $V_i$ the `vertical
arcs'
 and the  $H_i$ the
`horizontal' ones.

\begin{lemma}
\label{lemma:brokenarc}
Consider a broken arc in hyperbolic plane with endpoints $R,R'$ and 
with side lengths $s_1,d_1,\dots$,$s_r$,$d_r,s_{r+1}$. For any $D>0$,
there
exists a constant $K=K(D,r)$, depending only on $D$ and the number $r$
of horizontal arcs, so that, if $d_j>D$ for all $j$, we have
$$d(R,R')>\sum d_j +\sum s_j - K.$$
If $D'>D$, then $K(D',r)<K(D,r)$. 
\end{lemma}

 In the proof we use the following facts about universal constants for 
hyperbolic triangles, which can be deduced from the property that
hyperbolic triangles are thin, see for example~\cite{GdH}. \\
\noindent { (I)} There exists a positive constant $K(\theta_0)$ so that
for any hyperbolic triangle with side lengths $a,b,c$ and angle $\theta$
opposite to $c$ satisfying $\theta\geq \theta_0>0$, we have
$c>a+b-K(\theta_0).$ Moreover, if $\th'_0>\th_0$, then
$K(\th'_0)<K(\th_0)$.\\
\noindent { (II)} Given $D>0$, there exists a constant
$\theta_0=\theta_0(D)$ so that for any hyperbolic triangle with one side
of length $d\geq D$ and angles $\pi/2, \theta$ on this side, we have
$\theta\leq \theta_0$. If $D'>D$, then $\theta_0(D')<\theta_0(D)$.

\begin{proof}[Proof of Lemma~\ref{lemma:brokenarc}]
Consider $D>0$. The proof will be by induction on $r$.
For $r=1$ we have a broken arc with three arcs $V_1,H_1,V_2$ with
lengths $s_1,d_1,s_2$; denote by $Q,Q'$ the vertices of the arc $H_1$.
Since $d_1 > D$, by (II), there exists $\theta_0$ so that the angle
$QQ'R$ is less than  $\theta_0$; therefore the angle $RQ'R'$ is
greater than $\theta_1=\pi/2-\theta_0$. Applying (I) to the triangles
$RQQ'$ and $RQ'R'$, we have
$$d(R,R')>d(R,Q')+d(Q',R')-K(\theta_1)>s_1+d_1+s_2-K(\pi/2)-K(\theta_1)$$
so that we can take $K(D,1)=K(\pi/2)+K(\theta_1)$.

Now consider a broken arc with arc lengths
$s_1,d_1,\dots,s_r,d_r,s_{r+1}$ and  $d_i>D$. Denote by $Q,Q'$ the
vertices of the arc $d_r$. Since $d_r\geq D$, the angle $R'QQ'$ is
smaller than $\theta_0$. Since $R,R'$ are on different sides of the line
containing the vertical  segment $V_r$, the angle $RQR'$ is greater than
$\theta_1=\pi/2-\theta_0$. Applying (I) to the triangles $RQR'$
and $QQ'R'$ and using the induction hypothesis we get
$$d(R,R')>\sum_{j=1}^r d_j+\sum_{j=1}^{r+1} s_j
-K(D,r-1)-K(\pi/2)-K(\th_1).$$
So we can take $K(D,r)=r\left( K(\frac{\pi}{2})+K(\theta_1)\right)$.  If
$D'>D$,  by (I) and (II), we have that $K(D',r)<K(D,r)$.
\end{proof}

 Now let $\g$ be a closed geodesic on a hyperbolic surface $\rho$, and
let the geodesics $\{\a_i\}$ be a pants decomposition. We shall use  the
$\{\a_i\}$ to construct a broken arc
$BA_{\g}(\rho)$ associated to $\g$, as illustrated in Figure 2.
\begin{figure}[ht!]\small
\psfrag{C_1}{$\tilde C_1$}
\psfrag{C_r+1}{$\tilde C_r+1$}
\psfrag{C_2}{$\tilde C_2$}
\psfrag{C_r+1}{$\tilde C_{r+1}$}
\psfrag{H_1}{$H_1$}
\psfrag{H_2}{$H_2$}
\psfrag{H_r}{$H_r$}
\psfrag{H_r+1}{$H_{r+1}$}
\psfrag{Q_1^+}{$Q_1^+$}
\psfrag{Q_1^-}{$Q_1^-$}
\psfrag{Q_0^+}{$Q_0^+$}
\psfrag{r_1}{$\}r_1$}
\psfrag{Q_r^+}{$Q_r^+$}
\psfrag{Q_i^+}{$Q_i^+$}
\psfrag{Q_{i-1}^+}{$Q_{i-1}^+$}
\psfrag{Q_{i+1}^-}{$Q_{i+1}^-$}
\psfrag{V_2}{$V_2$}
\psfrag{V_{i+1}}{$V_{i+1}$}
\psfrag{R}{$\tilde Q$}
\psfrag{R'}{$\tilde Q'$}
\centerline{\includegraphics[height=6cm,width=12cm]{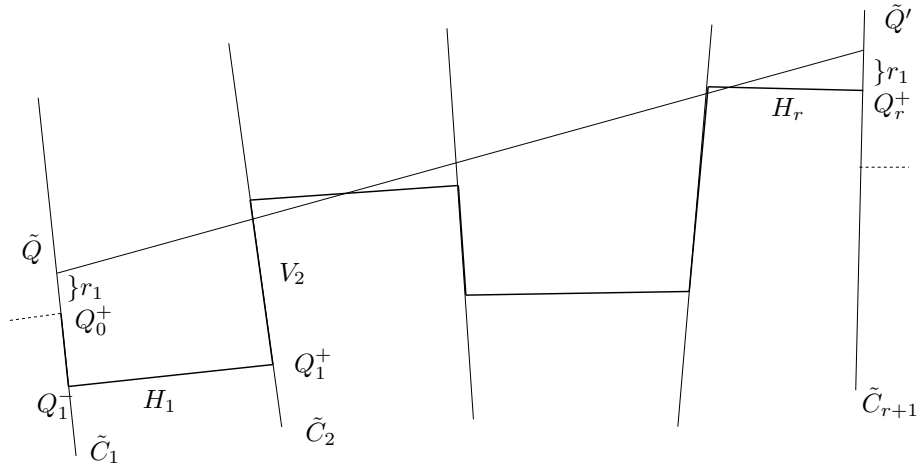}}
\caption{Broken arc with $r=4$ }
\end{figure}
  Fix  an orientation on $\g$ and let
$Q$ be  an intersection point of $\g$ with a pants curve. Let  $\tilde
\g$ be the lift of $\g$ through a lift $\tilde Q$ of $Q$.  Let
$\tilde{C_1},\dots,\tilde C_{r+1}$ be the lifts of the geodesics
$\{\a_i\}$
which are intersected, in order, by $\tilde{\g}$, so that $\tilde
C_1\cap \tilde\g=\tilde Q$ and  $\tilde C_{r+1}$ is the image of $\tilde
C_1$ under the covering translation $\tau$ corresponding to $\g$. Thus,
if we
denote by $\tilde Q'$ the intersection of $\tilde \g$  with $\tilde
C_{r+1}$,  the geodesic segment $\tilde Q \tilde Q'$ projects   onto
$\g$.   For $i=1,\dots,r$, consider the common perpendicular segment to
$\tilde{C_i},
\tilde C_{i+1}$, with endpoints denoted by $Q_i^-,Q_i^+$; and finally,
let $Q_0^+=\tau^{-1}
(Q_r^+)$. Then we define $BA_{\g}(\rho) $ to be the 
broken arc with vertical arcs the segments
$Q_0^+ Q_1^-, Q_1^+Q_2^-, \dots, Q_{r-1}^+Q_r^-$, and  horizontal arcs
the
segments $Q_1^-Q_1^+, \dots, Q_r^-Q_r^+$. Denote by $s_i$ the lengths of
the vertical
arcs and by $d_i$ the
lengths of the horizontal arcs.   
  The horizontal segments project onto geodesic
segments which are the perpendiculars either between two boundary
components, or from one boundary component to itself, of one of the
pairs of pants.  Their length will be studied in the next subsection.
The vertical segments  project onto arcs contained in the pants
geodesics $\a_j$. If the segment $Q_i^+Q_{i+1}^-$  projects, say, onto 
$\a_1$, then its length  is of the form
\begin{equation}
\label{eq:lengthVA}
s_i=|n_i l_{\a_1} +t_{\a_1}+e_i|
\end{equation}
where $n_i\in  \ZZ$ depends on the combinatorics of $\g$ relative to the
pants decomposition (related to how many times $\g$ wraps around
$\a_1$), 
 and 
$e_i$ is a number smaller in absolute value than $l_{\a_1}$ which  
depends on the combinatorics of $\g$ and on the geometry of the two
pairs of pants meeting along $\a_1$.  For our purposes we will not need
more details about $e_i$, see \cite{SW} for more explanation.   

We remark that the endpoints of $BA_{\g}(\rho)$ do not necessarily
coincide with those of $\tilde\g$, but we can consider another broken
arc $\overline{ BA}_{\g}$ with the same endpoints as $\tilde\g$ by just
changing the first vertical segment $Q_0^+Q_1^-$ to $\tilde QQ_1^-$ and
adding at the end the  vertical segment $Q_r^+\tilde Q'$.  To control
the lengths of these two new segments, we use the following lemma.

\begin{lemma}
\label{lemma:firstVA}
With the above notation, suppose that $\tilde C_i$ projects onto $\a_k$
and denote $\tilde Q_i=\tilde \g \cap \tilde C_i$. Then, either $\tilde
Q_i$ is between $Q_{i-1}^+$ and $Q_{i}^-$,  or the minimum of the
distances  $ d(Q_{i-1}^+,\tilde Q_i),   d(Q_{i}^-,\tilde Q_i)$ is less
than $l_{\a_k}$.  
\end{lemma}
\begin{proof}
Suppose that $\tilde Q_i$ is not between $Q_{i-1}^+$ and $Q_{i}^-$ and
that both distances $   d(Q_{i-1}^+,\tilde Q_i),  d(Q_{i}^-,\tilde Q_i)$
are greater than $l_{\a_k}$.  Then, applying the covering transformation
corresponding to $\a_k$ (or to $\a_k^{-1}$) to the segments
$H_{i-1},H_i$, we obtain two new  segments $H'_{i-1},H'_i$ which are
closer to $\tilde\g$ than $H_{i-1},H_i$. The lines containing these
segments are disjoint from $\tilde C_{i-1}, \tilde C_{i+1}$ respectively
and hence necessarily both intersect $\tilde \g$. Therefore they
determine, together with    the lines $\tilde \g, \tilde C_i$,   two
right-angled triangles.  One of them has  angle sum greater than $\pi$,
so we have a contradiction. (There is a similar argument when one of the
distances $   d(Q_{i-1}^+,\tilde Q_i),  d(Q_{i}^-,\tilde Q_i)$ is equal
to $l_{\a_k}$.)
\end{proof}
As a consequence,  if $\tilde C_1$ projects over $\a_{i_1}$ and  if
$\bar s_1,\bar s_{r+1}$ are the lengths of the first and last vertical
segments of  $\overline{ BA}_{\g}$, then    $\bar s_1+\bar s_{r+1}$ is
either equal to $s_1$ or to $s_1+2 r_1$, where $r_1<l_{\a_{i_1}}$. Thus,
by Lemmas~\ref{lemma:brokenarc} and \ref{lemma:firstVA}, we can
approximate the length of $\g$ by the length of the broken arc $BA_{\g}$
within an error of $K(D,r)+2l_{\a_{i_1}}$.  
 
\rk{Remark} A straightforward generalisation of the above
construction allows one to associate
a broken arc with any (not necessarily closed) geodesic, and also with a
geodesic arc with endpoints on the pants curves. (For a geodesic arc, to
determine the first vertical arc, prolong the geodesic  in the negative
direction until it crosses 
the next pants curve.) Then, we can use Lemmas~\ref{lemma:brokenarc} and
\ref{lemma:firstVA} to estimate this length from the length of this
broken arc. This is useful when $\nu$ is irrational, see Section 5.4.

\subsection{Geometry of a pair of pants}

We now estimate the lengths of the common
perpendicular segments between two curves of the pants decomposition. In
the situation to be considered, these segments will be sufficiently long 
to apply Lemma~\ref{lemma:brokenarc}.

It is useful to refine slightly the notation $f(s)\sim g(s)$ as $ s \to
s_0$ defined on p.~\pageref{page:order}. For $f,g$  real valued
functions we write  $f\approx g$ to mean that $\lim_{s \to s_0} f/g$
exists and is strictly positive.
Clearly, $f\approx g $ is slightly stronger than $f\sim g$.
However even if the limit  does not exist, if $f\sim g$, and if both
functions tend to  either $ 0$ or  $\infty$,  then  $ \lim_{s \to s_0}
\log f/\log g$ does exist and equals  $1$. This fact is crucial for our
results. We collect this and other elementary
properties in the next lemma.  We also recall the notation $f=O(g)$ as
$s\to s_0$ meaning that $f/g$ is bounded when $s\to s_0$, and  $f=o(g)$
as $s\to s_0$ meaning that $f/g \to 0$  when $s\to s_0$.

\begin{lemma}
\label{lemma:order}\begin{itemize}
\item[\rm (a)]  $ f \sim g $ is equivalent to  $\log f =\log g + O(1)$. 
\item[\rm (b)] If $f,g \to \infty$ or $0$ and $f\sim g$, then $\lim
(\log
f/\log g) =1$.
\item[\rm (c)]  $f\approx g$ is equivalent to $f=a g +o(g)$, with $a>0$.
\end{itemize}\end{lemma}
\begin{proof} (a)\qua If there exists  $0<k<K$ with $k<f/g<K$, then
taking logarithms we get that  
$
\log k <\log f -\log g <\log K.
$
The converse is also clear   by exponentiating $\log f =\log g +
O(1)$.

(b)\qua Since $g\to 0$ or $\infty$, then $\log g \to -\infty$ or $+\infty$
respectively, and in both cases $O(1)/\log g \to 0$. Then,  dividing 
$\log f =\log g + O(1)$ by $\log g$, we get the result. 
Part (c) is immediate from the definitions. 
\end{proof}

\begin{lemma}
\label{lemma:pairofpants} 
Consider a pair of pants $P$ with boundary components $B_1,B_2$, $B_3$ 
of lengths $l_1,l_2,l_3$. For any  $i,j \in \{1,2,3\}$, let $H_{ij}$ be
the common perpendicular arc to the boundary components $B_i,B_j$, with
length $d_{ij}$. Suppose that each of $l_1,l_2,l_3$ either tends to zero
or is bounded above. Then,   for any $i,j$, we have
$$
d_{ij} =\log \frac{1}{l_i} +\log\frac{1}{l_j} +O(1).
$$ 
 \end{lemma}
\begin{proof}
The pair of pants $P$ is made up by gluing two isometric right angle
hexagons with alternate sides of  lengths $l_1/2,l_2/2,l_3/2$. For
$i\not=j$, the segments $H_{ij}$ are the remaining sides. The segment
$H_{ii}$ is the union of the common  perpendicular segments in the two
hexagons between the side contained in $B_i$ and its  opposite side.
We therefore obtain the trigonometric formulae:
$$
\ch d_{ij}=\frac{\ch \frac{l_k}{2}+\ch \frac{l_i}{2}\ch \frac{l_j}{2}
}{\sh \frac{l_i}{2}\sh \frac{l_j}{2}}, \quad \ch \frac{d_{ii}}{2}=\sh
d_{ij}\sh \frac{l_j}{2}.
 $$
For $i\not=j$  we deduce that $\ch d_{ij}\approx
\frac{1}{l_il_j}$ as $(l_1,l_2,l_3)\to (0,0,0)$. Thus there exits $a>0$
so that 
$$
\ch d_{ij} =\frac{a}{l_il_j}+o\left( \frac{1}{l_il_j} \right). 
$$ 
Since $d_{ij}\to \infty$, $e^{-d_{ij}}$ is bounded, and so 
 $$
e^{d_{ij}} =\frac{2 a}{l_il_j}+o\left( \frac{1}{l_il_j} \right). 
$$ 
The result follows from Lemma~\ref{lemma:order} (a).

For the case $i=j$, 
 we   have that  $\sh d_{ij}\approx \frac{1}{l_il_j}$ as
$(l_1,l_2,l_3)\to (0,0,0)$ (because $d_{ij}\to \infty$ and in that case
$\sh d_{ij}\approx \ch d_{ij}$). Then, from the above formula for
$\ch(d_{ii}/2)$, we have that  
$\ch(d_{ii}/2) \approx 1/l_i$. As before, $d_{ii}/2 =\log(1/l_i)
+O(1)$ and so we get the result. 

We can check that the same works when some or all the $l_i$ do not tend
to zero but are still bounded above.
\end{proof}

\subsection{ Proof of Propositions~\ref{prop:geom.prop.LM} and
\ref{prop:lengthestimate} }
\label{sec:proofProps}

Proposition~ \ref{prop:lengthestimate} is used in the proof of
Proposition~\ref{prop:geom.prop.LM} (c), so we follow this order below.
We remark that in our proofs  the hypothesis of being in the line of
minima is always used 
in the same way: simply compare the value of $F_s$ at its minimum $m_s$
and at some other point.

\begin{proof}[Proof of Proposition~\ref{prop:geom.prop.LM}(a)]
Consider $i$ with $a_i\not=0$. Given $\epsilon >0$, consider a
hyperbolic surface $\rho_{\epsilon}$ so that the length of any $\a_i$ is
equal to $(a_i/4)\epsilon$. Take $s_0={\rm
min}\{\frac{1}{2},\frac{a_i\epsilon}{4l_{\nu}(\rho_{\epsilon})}\}$. 
Then, for  $s<s_0$, we have
\[
\begin{array}{ccl}
F_s(\rho_{\epsilon})&=&(1-s)\sum_{j=1}^N a_j l_{\a_j}(\rho_{\epsilon})+s
l_{\nu }(\rho_{\epsilon}) \\&=&
(1-s)(\sum_{j=1}^Na_j)\frac{a_i}{4}\epsilon+s l_{\nu}(\rho_{\epsilon})
<\frac{a_i}{4}\epsilon+\frac{a_i}{4}\epsilon=\frac{a_i}{2}\epsilon.
\end{array}
\]
Hence, for the minimum point $m_s$ of $F_s$, we have
$$(1-s)a_i l_{\a_i}(m_s)\leq F_s(m_s)<(a_i/2)\epsilon.$$
Since $1-s>1/2$, we have that $l_{\a_i}(m_s)<\epsilon$. Since, by hypothesis, $a_i\not=0$ for all $i$, we have the result.
\end{proof}

\begin{proof}[Proof of Proposition~\ref{prop:geom.prop.LM}(b)]  
We shall prove that all the twists about $\a_r$  are bounded when $s\to
0$. Suppose not; renumbering, we may assume that   $|t_{\a_1}(m_s)|$  is 
not bounded. Suppose moreover that, up to subsequence, $t_{\a_1}\to +
\infty$. (The proof is the same if $t_{\a_1} \to -\infty$.) 

For each $s$, consider the point $\rho_s$ obtained from $m_s$ by  
twisting (earthquaking) by $-t_{\a_1}(m_s)$  about $\a_1$. This new
surface has the same   Fenchel-Nielsen
coordinates with respect to $\{\a_i\}$ as $m_s$ (for a fixed choice of
dual curves) except that  $t_{\a_1}(\rho_s)=0$. We shall
prove that, if $s$ is small enough, then $l_{\nu}(m_s)>l_{\nu}(\rho_s)$.
Since the lengths of $\a_i$ are the same at both points, we will have
that $F_s(m_s)>F_s(\rho_s)$, which is a contradiction.

We make the following argument for all the curves $\b$ in the support of
 $\nu$.  
Since the lengths of all the pants curves tend to zero, by
Lemma~\ref{lemma:pairofpants} there exists
$s_0>0$ so that for all $s<s_0$, all the horizontal arc lengths $d_{ij}$
of the broken
arc $BA_{\b}$ are bounded below by  some constant $D$, and therefore by
Lemma~\ref{lemma:brokenarc} (and Lemma~\ref{lemma:firstVA}), there
exists a constant $K$ so that 
$l_{BA_{\b}}(m_s)-l_{\b}(m_s)<K$ 
  for all $s<s_0$. 
We will prove that
 \begin{equation}
\label{eq:comparisonBA}
 l_{BA_{\b}}(m_s)-l_{BA_{\b}}(\rho_s) \to \infty  \; \ {\rm as} \; \ 
s\to 0.
\end{equation}
Assuming (\ref{eq:comparisonBA}), there  exists $0<s_1<s_0$ so that for
all $s<s_1$, $l_{BA_{\b}}(m_s)-l_{BA_{\b}}(\rho_s)>K$, and then
$$l_{\b}(m_s)>l_{BA_{\b}}(m_s)-K> l_{BA_{\b}}(\rho_s)
>l_{\b}(\rho_s).$$
Summing over all  curves $\b_i$, we obtain
$l_{\nu}(m_s)>l_{\nu}(\rho_s)$.

It is left to prove (\ref{eq:comparisonBA}). 
We compare the broken arcs $BA_{\b}$ at the points $m_s$ and
$\rho_s$. The horizontal arcs have the same length in both broken arcs.
In fact only the vertical arcs projecting over the curve $\a_1$ change 
length. There are  $i(\b,\a_1)$ of such segments,  with lengths
$s_1,\dots, s_{i(\b,\a_1)}$ at $m_s$, where 
$ s_j=| n_j l_{\a_1} +t_{\a_1} +e_j|$. On the other hand, the lengths  
of these segments  at $\rho_s$ are   
$s_j'=|n_j l_{\a_1} +e_j|.$ (We remark that $n_1+\dots+n_{i(\b,\a_1)}$
is called the {\it wrapping number} of $\b$ around $\a_1$.) 
Then 
$$l_{BA_{\b}}(m_s)-l_{BA_{\b}}(\rho_s)= \sum (s_j-s_j')\geq
i(\b,\a_1) t_{\a_1} -2\sum_{j=1}^{i(\b,\a_1)}|n_jl_{\a_1} +e_j| .$$
Since $|n_jl_{\a_1} +e_j|$ is bounded for all $j$, this expression tends
to
infinity as $s\to 0$ as required.
\end{proof}

\begin{proof}[Proof of Proposition~\ref{prop:lengthestimate}]
Consider the broken arc $BA_{\g}(\rho_n)$ associated to $\g$. Since all
the lengths $l_{\a_i}$ tend to 0 or are bounded above when $n\to \infty$,
then, for $n$ big enough,  all horizontal arcs of this broken arc are
greater than some given $D$. Then we can  use
Lemmas~\ref{lemma:brokenarc} and \ref{lemma:firstVA}
to estimate the length of $\g$ and we have 
$
\sum d_i +\sum s_i -K < l_{\g}(\rho_n) < \sum d_i +\sum s_i.
$
Since  all the twists $t_{\a_i}$ are bounded, the lengths of the
vertical arcs are bounded. On the other hand, using Lemma \ref{lemma:pairofpants} to estimate the lengths $d_i$ of the horizontal arcs, and collecting terms together, we have   
$$
  l_{\g}(\rho_n)=\sum d_i
+O(1)= \sum_j r_j\log \frac{1}{l_{\a_j}(\rho_n)} +O(1),
 $$ 
where $r_j$ is the number of times that the projections of the horizontal arcs end in the geodesic $\a_j$.  Then this number is  equal to $2i(\a_j,\g)$, and we get the desired result.
\end{proof}

\begin{proof}[Proof of Proposition~\ref{prop:geom.prop.LM}(c)]
 We need to prove that, on the line of minima ${\cal L}_{\mu,\nu}$, all
the lengths $l_{\a_i}$ have  the same order as $s$, when $s\to 0$. That
is,  there exists $\epsilon >0$
 and   positive constants $k_i<K_i$
for any $i$, so that 
$$k_i <\frac{l_{\a_i}(m_s)}{s}<K_i $$
for all $s<\epsilon$. Suppose not; then, up to subsequence,  there
exists   some curve $\a_j$ so that, when $s\to 0$, then
$l_{\a_j}(m_{s})/s$ tends either to 0 or to $\infty$. 

We construct a new sequence of surfaces $\rho_{s}$ defined by
the Fenchel-Nielsen coordinates  $l_{\a_i}=s$ and $t_{\a_i}=0$, for all
$i$.  
We compare $F_{s}(m_{s})$ with $F_{s}(\rho_{s})$, and show that 
$F_{s}(m_{s})-F_{s}(\rho_{s}) >0$, which will be a contradiction because
$m_s$ is the minimum of the function $F_s$.

By Proposition~\ref{prop:geom.prop.LM} (a) and (b), along $m_{s}$ the
lengths $l_{\a_i}$ tend to 0 and the twists  $t_{\a_i}$ are bounded.
Applying Proposition~\ref{prop:lengthestimate} to estimate the
length of the curves $\b_i$ we find
$$
\textstyle{F_s(m_s)=(1-s)\sum_i a_i l_{\a_i}(m_{s}) +s\sum_i b_i \left(
2\sum_j
i(\a_j,\b_i) \log\frac{1}{l_{\a_j}(m_s)}  +O(1) \right).}
$$
On the other hand, the sequence $\rho_s$ also satisfies the hypothesis
of Proposition~\ref{prop:lengthestimate}, and so we can also use this
proposition to estimate the length of $\b_i$, giving 
$$
\textstyle{F_s(\rho_s)=(1-s)\sum_i a_i s +s\sum_i b_i \left( 2\sum_j
i(\a_j,\b_i)
\log\frac{1}{s}  +O(1) \right).}
$$
Then
\begin{eqnarray*} 
\lefteqn{\frac{F_s(m_s)-F_s(\rho_s)}{s}}
 \hspace{-.2cm}\\
&&=(1-s)\textstyle{\sum_i a_i \frac{l_{\a_i}(m_s)-s}{s}}\\
 && \textstyle{ \hspace{1.5cm}{+\sum_ib_i\left( 2\sum_j
i(\a_j,\b_i)\left( \log
\frac{1}{l_{\a_j}(m_s)} -\log \frac{1}{s} \right) +O(1) \right)}}\\ 
&&\textstyle{= (1-s)\sum_i a_i \frac{l_{\a_i}(m_s)}{s} +\sum_ib_i\left(
2\sum_j
i(\a_j,\b_i)\left( \log \frac{s}{l_{\a_j}(m_s)}   \right)  \right)
+O(1)}\\
&&\textstyle{ =\sum_i \left( (1-s) a_i \frac{l_{\a_i}(m_s)}{s} +C_i \log
\frac{s}{l_{\a_i}(m_s)}  \right) +O(1),} 
\end{eqnarray*}
where in the last equality we have rearranged the second group of
summands,
and  $C_i$ are some positive coefficients.

Now, if there is some $i$ so that $l_{\a_i}(m_s)/s \to 0$, then $\log
(s/l_{\a_i}(m_s) )\to \infty$. On the other hand, if there is some $i$
so
that 
$l_{\a_i}(m_s)/s \to \infty $ then $\log (s/l_{\a_i}(m_s) )\to -\infty$,
but any positive linear combination of $l_{\a_i}(m_s)/s$ and  $\log
(s/l_{\a_i}(m_s))$ tends to infinity. 
 
Hence, $(F_s(m_s)-F_s(\rho_s))/s$ tends to $+\infty$ and therefore
$F_s(m_s)-F_s(\rho_s)>0$ for sufficiently small $s$, obtaining the
desired
contradiction.
\end{proof}

\subsection{The proof of  Theorem~\ref{thm:linetobarycenter}  when  
$\nu$ is irrational}
\label{sec:irrational}  

Finally we discuss the proof of   Theorem~\ref{thm:linetobarycenter} 
when   $\nu$ is irrational. All that is needed is to extend
Propositions~\ref{prop:geom.prop.LM}
and~\ref{prop:lengthestimate} 
to the case in which $\nu$ and $\gamma$ respectively are general
measured laminations. 

First, consider the effect of replacing $\g$ in
Proposition~\ref{prop:lengthestimate} 
by an irrational lamination $\eta$. There exists a sequence of rational
laminations  $c_k\g_k$ converging to $\eta$, with $\g_k$ simple closed
curves and $c_k\to 0$. Then $l_{\eta}(\rho_n)=\lim_{k\to \infty}
c_kl_{\g_k}(\rho_n)$. We can compute this limit by using the expression
obtained in Proposition~\ref{prop:lengthestimate} for closed curves.
Since $\lim_{k\to\infty}i(c_k\g_k,\a)=i(\eta,\a)$, 
then we only need to check that the error  (which is a  function 
$f_k(n)$, bounded for fixed $k$ as $n\to \infty$) stays bounded when
$k\to \infty$.  By careful inspection of  this error (in the  proofs of
Lemmas~\ref{lemma:brokenarc}, \ref{lemma:pairofpants} and
Proposition~\ref{prop:lengthestimate}),  we see that it depends linearly
on the intersection number of $\g_k$  with $\sum\a_i$ and  on the
wrapping numbers of $\g_k$ around $\a_i$. In both cases, these numbers,
after scaling with $c_k$,   converge when $k\to \infty$ (to the
intersection number of $\eta$ with $\sum \a_i$ and to the {\it twisting
numbers} of $\eta$ around $\a_i$ respectively). This proves
Proposition~\ref{prop:lengthestimate}.

The proof of Proposition~\ref{prop:geom.prop.LM} (b) for $\nu$
irrational, uses the same kind of arguments. The proof of part (c)  is 
unchanged, once we have the stronger version of
Proposition~\ref{prop:lengthestimate}.

Alternatively,  we can compute the length of an irrational lamination 
$\eta$ from its definition (see \cite{KerckEA}), as the  integral over
the surface of the product measure $d\eta \times d l$, where $dl$ is the
length measure along the leaves of $\eta$. We can cover the surface with
thin rectangles, so that the length of arcs of $\eta$ intersecting one
rectangle  are almost equal, and we approximate the length of one of
these geodesic arcs by  the length of a broken arc (notice that we need
the remark after Lemma~\ref{lemma:firstVA} to do this).    In this way 
we can prove both Propositions~\ref{prop:lengthestimate} and
\ref{prop:geom.prop.LM} (b).  
\qed

\section{Non-pants decomposition case}
\label{sec:npdc}
 We now investigate the modifications needed to the above work if
  $\a_1,\dots , \a_N$ is not a pants decomposition.
 The problem is that we no longer have full control over the geometry of
the complement in $S$ of the  curve system $\a_1,\dots , \a_N$. The
hyperbolic structures on at least some  components of the complement
might themselves diverge, in other words, the estimate of
Proposition~\ref{prop:lengthestimate} for the lengths of arbitrary
closed geodesics may no longer hold. Without ruling out this
possibility, we show that the divergences in question must be of a lower
order than those caused by the shrinking of the curves $\a_i$.
A precise statement is made in Corollary~\ref{cor:estimateNP}.

In more detail, we proceed as follows.  With minimal changes we still
  can prove that $l_{\a_i}\to 0$, with the same order as $s$, and that
  the twists $t_{\a_i}$ are bounded
  (Proposition~\ref{prop:gpLM2}). This implies that, if there exists a
  limiting lamination $\eta$, then its support is either disjoint from
  or contains the curves $\a_i$. Next we prove that there is in fact
  no other lamination contained in the limiting lamination, and
  therefore $[\eta]=[a'_1\a_1+\dots+a'_N\a_N]$, for some coefficients
  $a'_i\geq 0$ (Proposition~\ref{prop:support}). To compute these
  coefficients we need to compare the lengths of two closed
  geodesics. Even though we no longer have
  Proposition~\ref{prop:lengthestimate}, we can still extend the
  curves $\a_i$ to a pants decomposition and estimate the length of
  the dual curves, (Proposition~\ref{prop:estimatedual}). This is
  enough to compute the coefficients $a'_i$ and prove
  Theorem~\ref{thm:linetobarycenter}.  A posteriori we obtain an
  estimate for the length of closed geodesics along the line of minima
  in Corollary~\ref{cor:estimateNP}.

\begin{prop}
\label{prop:gpLM2}
Suppose that $\mu=\sum a_i\a_i$, $\nu=\sum b_i\b_i$ are two measured
laminations
(where $\{\a_1,\dots,\a_N\}$ is  not necessarily a pants decomposition),
and that $a_i>0$
for all $i$. Let $m_s$ be the minimum point of the function $F_s$. Then: 
\begin{itemize}
\item[\rm(a)] for any $i$, $\lim_{s\to 0}l_{\a_i}(m_s)=0$;
\item[\rm(b)] for any $i$, $|t_{\a_i}(m_s)|$ is bounded when $s\to 0$;
\item[\rm(c)] for all $i,j$, $l_{\a_i}(m_s)\sim l_{\a_j}(m_s)\sim s$ as
$s\to 0$.
\end{itemize}\end{prop}


\rk{Remark}
The twist parameter $t_{\a_i}$ is a real
parameter which
determines how the surface is glued along $\a_i$. It is determined
up to the choice of an initial surface on which the twist is zero.
Changing the initial surface results in  an additive change to the twist
parameter, so (b) above is independent of this choice.

\begin{proof}
The proof of (a) is exactly the same as that of
Proposition~\ref{prop:geom.prop.LM}(a).

For (b), consider a broken arc associated to the curve $\b_i$
relative to the curves $\a_1,\dots, \a_N$. Even if these curves are not
a pants decomposition, the definition of broken
arc given in Section~\ref{sec:brokenarcs} makes sense. The horizontal
segments 
now project onto arcs perpendicular to two of the $\a_i$. If the length
of each $\a_i$ is sufficiently small, the horizontal segments are
greater
that some given positive constant $D$ and we can apply
Lemmas~\ref{lemma:brokenarc} and \ref{lemma:firstVA} to approximate the
length of $\b_i$ with the
length of the corresponding  broken arc. The length of the vertical arc
projecting
over $\a_i$  is $|  t_{\a_i}+ r|$,  where $r$  depends only
on the curve $\b_i $ and the system $\a_1,\dots,\a_N$. (To see this,
think about obtaining the given surface from an initial surface with 
$ t_{\a_i} =0$  by twisting about $\a_i$.) We can
therefore argue exactly as in the proof of
Proposition~\ref{prop:geom.prop.LM} (b):
consider points $\rho_s$ obtained from the points $m_s$ on the line of
minima  by twisting by $-t_{\a_i}(m_s)$ about each $\a_i$. As
in that proposition, we obtain that $F_s(m_s)>F_s(\rho_s)$ for small
enough $s$.

For (c) we follow the same  argument as in the proof of
Proposition~\ref{prop:geom.prop.LM} (c). Let $\g_1,\dots,\g_K$ be simple
closed curves so that $\a_i,\g_j$ are a pants decomposition. Consider 
the surfaces $\rho_s$ defined by $l_{\a_i}(\rho_s)=s$,
$l_{\g_i}(\rho_s)=1$, $t_{\a_i}(\rho_s)=t_{\g_i}(\rho_s)=0$. The family 
$\rho_s$ satisfies the hypothesis of
Proposition~\ref{prop:lengthestimate}, and so we can estimate the length
of $\b_i$ as
  $$
l_{\b_i}(\rho_s)= 2\sum_j i(\a_j,\b_i) \log\frac{1}{s}  +O(1) .
$$

Now each curve $\a_i$ is contained in an embedded annular collar of
width at least $ 2 \log (1/l_{\a_i})$. Using the contribution of the
these collars gives the following rough estimation for the length
of $\b_i$ at $m_s$:
$$
l_{\b_i}(m_s)= 2\sum_j i(\a_j,\b_i) \log\frac{1}{l_{\a_j}(m_s)}  +f(s) ,
$$
where $f(s)$ is a positive function that might tend to infinity. 
Now, as in the proof of  Proposition~\ref{prop:geom.prop.LM} (c), we
have
$$
\textstyle{\frac{F_s(m_s)-F_s(\rho_s)}{s}=
\sum_i \left( (1-s) a_i \frac{l_{\a_i}(m_s)}{s} +C_i \log
\frac{s}{l_{\a_i}(m_s)}  \right)+f(s) -O(1)}
$$
where $f(s)-O(1)$ may be negative but is nevertheless bounded below. The
conclusion follows as in that proposition since, if $l_{\a_i}$ does not
have the same
order as $s$, the above group of summands always tends to $+\infty$.
\end{proof}

\begin{prop}
\label{prop:support}
Let $\mu$, $\nu$ and $m_s$ be as in {\rm Proposition~\ref{prop:gpLM2}}
and suppose that $m_s$ converges to a projective measured lamination
$[\eta]$. Then the support of $[\eta]$ is contained in the union of
$\a_1, \dots, \a_N$.
\end{prop}
\begin{proof} If the conclusion is false, then using
Proposition~\ref{prop:gpLM2} (a) we must have  $|\eta|\subset \a_1\cup
\dots\cup \a_N\cup |\d|$,
where $\d$ is a measured  lamination whose support is  disjoint from
the $\a_i$.  Since $\mu$ and $\nu$ fill up the surface and
$i(\mu,\d)=0$, it follows that $i(\nu,\d)\not=0$, and therefore some
curve  $\b \subset |\nu|$ intersects
$|\d|$. Let $\bar \kappa$ be a geodesic arc contained in $\b$,  intersecting $\d$,  running
from  some $\a_i$ to some $\a_j$ (where possibly  $\a_i=\a_j$)  and not
intersecting any other $\a_l$.  We take  open collar neighbourhoods
$A_i$ of the curves $\a_i$, of width $2\log (1/l_{\a_i})$,
  and 
let $\kappa = \bar \kappa  - \bar \kappa \cap (A_i \cup A_j)$, so that
$\kappa$ is a geodesic segment with endpoints on the  relevant
components
$\tilde \a_i,
\tilde \a_j$  of the boundaries $\partial A_i$ and $\partial A_j$. 
Note that      the boundary curves
$\tilde \a_i, \tilde \a_j$ of the collars $A_i,A_j$ have length $O(1)$.

We are going to prove  that $l_{\kappa}$, the length of the geodesic
segment $\kappa$, tends to infinity by comparing to the length of a
simple closed curve or curves we call the {\it double} of $\kappa$ (or
$\bar \kappa$). 
If $i \neq j$ (or if $i=j$ but $\kappa$ meets both boundary components of the collar $A_i$), then the double is the simple closed curve $\tilde \b$ created
by  going around $\tilde \a_i$,
then parallel to $\kappa$, around $\tilde \a_j$, and back parallel to
$\kappa$.
In the case that
$\a_i=\a_j$ and $\kappa$ intersects only one boundary component of the collar $A_i$, then $\kappa$ splits this boundary component, $\tilde \a_i$, into two arcs 
$\tilde \a_i'$ and $\tilde \a_i''$; we create two simple closed
curves  
$\tilde \b'=  \kappa \cup \tilde \a_i'$ and $\tilde \b''=  \kappa\cup
\tilde \a_i''$, and designate  $\tilde \b= \tilde \b' \cup
\tilde \b''$  the   double. 
 In Lemma~\ref{lemma:doublingarc} below, we show that $i(\tilde\b,
\d)\not=0$, where in  the second case we define   $i(\tilde\b,
\d) = i(\tilde \b' , \d)+i(\tilde\b'' ,
\d)$.

Since $m_s\to [\sum a'_i\a_i+\d]$,
the length on the surface  $m_s$ of the geodesic(s) isotopic to 
$\tilde\b$ must tend
to infinity as $s\to 0$.  We claim that the length of the arc $\kappa$
at $m_s$ tends to infinity with $\tilde\b$.  If $\a_i\not=\a_j$ then
$ 2l_{\kappa} + l_{\tilde \a_i}+ l_{\tilde \a_j} > l_{\tilde \b} $.
The lengths of the boundary curves $\tilde \a_i, \tilde\a_j$ are bounded
above (and below); since $ l_{\tilde \b}\to
\infty$, this forces $l_{\kappa} \to \infty $. A similar proof works  if 
$\a_i=\a_j$. 

Finally, we use the hypothesis that $m_s$ is the minimum of $F_s$ to
arrive to a contradiction. The argument is similar to others used above. 
Let $\g_1,\dots,\g_K$ be simple closed curves extending
$\a_1,\dots,\a_N$ to a pants decomposition, and fix a set of dual
curves. For each $s$ let $\rho_s$ be
the surface whose Fenchel-Nielsen coordinates with respect to these
choices   are 
$$
l_{\a_i}(\rho_s)=l_{\a_i}(m_s), \ l_{\g_i}(\rho_s)=1, \
t_{\a_i}(\rho_s)=
t_{\g_i}(\rho_s)=0.
$$
The surfaces $\rho_s$ satisfy the hypothesis of
Proposition~\ref{prop:lengthestimate}, and therefore 
$
l_{\b_i}(\rho_s)=\sum_j2i(\b_i,\a_j)\log(1/l_{\a_j}) +O(1)
$. 
On the other hand, we have 
\[
l_{\b_i}(m_s)=\sum_j2i(\b_i,\a_j)\log(1/l_{\a_j}) +f(s),
\]
where $f(s)$ is a positive function which tends to infinity  for those
curves $\b_i$ intersecting $|\d|$, since, by the above argument, some
arcs of some $\b_i$ outside the collars $A_j$ tend to infinity. 
Thus  
$$\textstyle{F_s(m_s)-F_s(\rho_s)= s \sum_i b_i
\left(l_{\b_i}(m_s)-l_{\b_i}(\rho_s)\right)=s\sum_ib_i(f(s)-O(1))}
$$
which is positive for  small enough $s$.
\end{proof}

The following lemma   was used in the above proof. We provide a proof, although the result describes a well-known construction.
\begin{lemma}
\label{lemma:doublingarc}
Let $\rho$ 
be a hyperbolic surface,   let  $\a_1,\a_2$  be  two disjoint simple
closed
geodesics, and   $\bar \kappa$ be a geodesic arc from $\a_1$ to $\a_2$.
Let $\d$
be a geodesic intersecting $\bar \kappa$ and    $\tilde  \b$ be  the
`double' of  $\bar \kappa$, as constructed
 in the proof  of 
{\rm Proposition~\ref{prop:support}}. Then $i(\tilde\b,\d)\not=0$.
\end{lemma}
\begin{proof}
Suppose that  $\a_1\not=\a_2$; by the construction of the curve
$\tilde\b$,  the curves  $ \a_1, \a_2,\tilde\b$ bound a pair of pants,
made up of a thin strip around the arc $\kappa$ (the part of $\bar
\kappa$ outside the annuli $A_i$), and the sub-annuli of $A_1$ and $A_2$
with boundaries $\a_i$ and $\tilde \a_i$ for $i=1,2$. 
Correspondingly, there is a pair of pants $P$  in our hyperbolic surface
bounded by $\a_1,\a_2$ 
 and the geodesic representing  $ \tilde\b$.  Now,  $\bar \kappa$ is an
arc
contained in $P$  joining $\a_1$ to $\a_2$. The lamination $\d$
intersects  $P$ in arcs which do not meet $\a_1,\a_2$, so 
running from the geodesic representative of $\tilde\b$ to itself;
therefore   each of these arcs
intersects   $\bar \kappa$ once and $\tilde\b$ twice so that
$i(\tilde\b,\d)=2i(\bar \kappa,\d) >0$.

In the case that $\a_1=\a_2$, remember that the
`double'  of  $\bar \kappa$ is the union of the two simple closed curves
$\tilde
\b',\tilde\b''$ described above. We have to show that
$i(\tilde\b',\d)+i(\tilde\b'',\d)\not=0$. Arguing much as above,  we
have that $\a_1$ and  the geodesic  representatives of $\tilde\b'$ and
$\tilde \b''$ bound a pair of pants containing $\bar \kappa$, and $\bar
\kappa$ joins $\a_1$ to itself. Since $\d$ intersects this pair of pants
in geodesic  arcs not intersecting $\a_1$,  each such arc  intersects  
$\tilde\b'\cup\tilde\b''$ twice and  the result follows.
\end{proof} 

\medskip 
 Finally, we estimate the length along the line of minima of curves dual
to the $\a_i$. Suppose $\a_i,\g_j$ is a pants decomposition of the
surface $S$, and let $\d_i$ be the dual curves. If $\d_i$ is dual to $\a_i$,
then these two curves intersect either once or twice.  
If $i(\a_i,\d_i)=1$, then $\a_i$ is on the boundary of just one pair of pants
$P$  (two boundary components of $P$ are glued
together along $\a_i$). We denote the other boundary component of $P$ by
$\omega $.  If $i(\a_i,\d_i)=2$, then $\a_i$ is on the boundary of two
different pants $P, P'$; let $\omega_1,\omega_2, \omega'_1,\omega'_2$
the other boundary components of $P,P'$, respectively. To simplify notation, in the following proposition we drop the indices in $\a_i,\d_i$.

\begin{lemma}
\label{lemma:dual} With the above notation, let $\d$ be the dual curve 
 to $\a$.\begin{itemize}
\item[\rm (a)] Suppose that $i(\a,\d)=1$ and that  $\rho_n$ is a
sequence of surfaces so that 
$l_{\a}(\rho_n)\to 0$ and  $|t_{\a}(\rho_n)|$ is bounded. Then
$$   
l_{\d}(\rho_n)= 2i(\a,\d) \log(1/l_{\a}(\rho_n))+(1/2)
l_{\omega}(\rho_n)+O(1)
.$$
\item[\rm (b)] Suppose that $i(\a,\d)=2$ and suppose that  $\rho_n$ is a
sequence of surfaces so that 
$l_{\a}(\rho_n)\to 0$ and $|t_{\a}(\rho_n)|$ is bounded. Then
$$   
l_{\d}(\rho_n)= 2 i(\a,\d) \log(1/l_{\a}(\rho_n))+ l_{|\omega|
}(\rho_n)+
l_{|\omega'|}(\rho_n)+O(1),
$$
where we denote by $l_{|\omega| }(\rho_n)$ the maximum of
$l_{\omega_1}(\rho_n)$ and $l_{\omega_2 }(\rho_n)$, and similarly with
$l_{|\omega'| }(\rho_n)$.\end{itemize}
\end{lemma}
\begin{proof}
 In both cases, the length of $\d$ depends on the lengths of $\a$ and
of the neigbouring pants curves and on the twist about $\a$. One can
calculate $\d$ explicitly,  however it is easier to use the broken arc
Lemma~\ref{lemma:brokenarc} to simplify the estimates.

If $i(\a,\d)=1$, let  $d$ be the distance in $P$ between the
two boundary components $\a', \a''$ projecting over $\a$.
The dual curve $\d$ can be approximated by a broken arc which wraps
part-way round $\a'$ and then follows the common  perpendicular from
$\a'$ to $\a''$, and finally wraps part-way round $\a''$.
Since we are assuming that $l_{\a} \to 0$ and
$|t_{\a}(\rho_n)|$ is bounded,   Lemma~\ref{lemma:brokenarc} gives the
approximation 
$l_{\d} = d+O(1)$.

If $i(\a,\d)=2$, let
 $b, b'$ be the lengths
of  the common perpendicular arcs from $\a$ to itself in  $P$ and $P'$.
In this case  the approximating broken arc has five segments; three
vertical segments which each wrap part-way round arcs which project to
$\a$, and  two horizontal segments which are just the common
perpendiculars from 
$\a$ to itself in $P$ and $P'$. Thus in this case 
Lemma~\ref{lemma:brokenarc} gives the approximation  $l_{\d} =
b+b'+O(1)$.

The proof is completed by using the trigonometric formulae in the proof
of Lemma~\ref{lemma:pairofpants} to estimate $d,b$ and $b'$.
For (a),  if
$l_{\omega}$ is bounded above, then  $d=2\log
(1/l_{\a}(\rho_n))+O(1)$; while if $l_{\omega}\to \infty$, then $\ch d
\approx
e^{l_{\omega}/2}/l_{\a}^2$, so  the result still holds.

For case (b), note that $P$ is made up of two right-angled hexagons
with alternate sides of lengths
$l_{\a}/2,l_{\omega_1}/2,l_{\omega_2}/2$ and that $b/2$ is the
distance between $a$ and its opposite side. Let $d_1$ be the length of
the side between the sides of lengths $l_{\a}/2$ and
$l_{\omega_1}/2$. We claim that $b= 2\log(1/l_{\a})
+l_{|\omega|}+O(1)$, from which, combined with a similar estimate for
$b'$, part (b) follows.

Since $l_{\a}\to 0$, we have $d_1\to \infty$, and so $\ch
d_1 \approx \sh d_1$. Thus
$$
\ch \frac{b}{2}\approx
\frac{\ch\frac{l_{\omega_2}}{2}+\ch\frac{l_{\a}}{2}
\ch\frac{l_{\omega_1}}{2}}{\sh\frac{l_{\a}}{2}\sh\frac{l_{\omega_1}}{2}}\sh
\frac{l_{\omega_1}}{2}
\approx \frac{\ch\frac{l_{\omega_2}}{2}+\ch
\frac{l_{\omega_1}}{2}}{l_{\a}}.
$$ 
Since $b \to \infty$, we have $ \ch b/2 \approx e^{b/2}$, so by
Lemma~\ref{lemma:order},  
\begin{equation}
\label{eq:formulab1}
 b/2=\log ( \ch(l_{\omega_2}/2)+\ch (l_{\omega_1}/2)  )+\log
(1/l_{\a}) +O(1).
\end{equation}
Now, expressing $\ch(l_{\omega_2}/2)+\ch (l_{\omega_1}/2) $ as
${\rm max}\{ \ch(l_{\omega_2}/2),\ch (l_{\omega_1}/2)\}+{\rm min}
\{\ch(l_{\omega_2}/2),\ch (l_{\omega_1}/2)\}$, we  easily obtain that 
$$\log ( \ch(l_{\omega_2}/2)+\ch (l_{\omega_1}/2)  ) = l_{|\omega|}
/2+O(1).$$
Applying
this to (\ref{eq:formulab1}) gives the  claim. 
\end{proof}

\begin{prop}
\label{prop:estimatedual}
Let $\a_1,\dots, \a_N, \g_1,\dots, \g_K$ be a pants decomposition of $S$
and let  $\d_i$ be the dual curve to $\a_i$. Let $\rho_n$ be a sequence
so that 
$l_{\a_i}\to 0$, $|t_{\a_i}|$ is bounded and $\rho_n\to [a'_1a_1+\dots
+a'_N\a_N]$. Then, for any $j$ with $a'_j\not=0$ we have 
$$
l_{\d_j}(\rho_n)=2i(\a_j,\d_j)\log (1/l_{\a_j}) +o(\log
(1/l_{\a_{j}})).
$$ 
\end{prop}
\begin{proof} 
Suppose that $\a_j,\d_j$ intersect twice (the proof is similar if they
intersect once).   From  Lemma~\ref{lemma:dual} (b) we have
$$
l_{\d_j}=2i(\a_j,\d_j) 
\log(1/l_{\a_j})+l_{|\omega|}+l_{|\omega'|}+O(1).
$$
 Since $\rho_n\to [a'_1a_1+\dots +a'_N\a_N]$, we have 
$l_{\d_j}/c_n\to a'_j i(\a_j,\d_j)$, for some sequence $c_n\to \infty$.
If
$\omega_1 $ is  one of the curves $\a_i$, then
$l_{\omega_1}(\rho_n) \to 0$. Otherwise, $\omega_1$ is
one of the curves $\g_i$, so it is disjoint from the curves $\a_i$, in
which case 
$l_{\omega_1}(\rho_n)/c_n \to i(\omega_1, a'_1\a_1+\dots+a'_N\a_N)=0$. 
The same holds for $\omega_2,\omega'_1,\omega'_2$. Thus
$$
\lim_n\frac{l_{\d_j}(\rho_n)}{c_n}=\lim_n
\frac{2i(\a_j,\d_j)\log(1/l_{\a_j}(\rho_n))}{c_n}= a'_j i(\a_j,\d_j), 
$$
which implies that $c_n\approx \log (1/l_{\a_j}(\rho_n))$, and therefore
$\frac{l_{|\omega|}(\rho_n)+l_{|\omega'|}(\rho_n)+O(1)}{\log
(1/l_{\a_j})}\to
0$. Thus $l_{|\omega|}(\rho_n)+l_{|\omega'|}(\rho_n)+O(1)=o(\log
(1/l_{\a_j}))$, which completes the proof.
\end{proof}

\medskip
  We can now complete the proof of Theorem~\ref{thm:linetobarycenter}.
First, continuing with the assumption that $\nu$ is rational,
we follow the method used in Section 4 for the pants decomposition case. 
  Proposition~\ref{prop:support} shows that the limit of any convergent
subsequence of minima $m_s$  is a projective lamination
$[a'_1\a_1+\dots+a'_N\a_N]$ for some $a'_i\geq 0$. 
Proposition~\ref{prop:gpLM2} (c)
 implies that $  \lim \frac{\log(1/l_{\a_i})}{\log(1/l_{\a_j})}=1$, for
all $i,j$ and therefore, using Proposition~\ref{prop:estimatedual},
we argue as in Section 4 to get that 
$\lim_{s\to0} \frac{l_{\d_i}(m_s)}{l_{\d_j}(m_s)}= i(\sum_k
\a_k,\d_i)/i(\sum_k \a_k,\d_j )$, so that $a'_k=1$ for all $k$.
Thus the limit is independent of the subsequence, and the result follows
by compactness of  $\teich(S) \cup \PML$. 
Finally, to complete the proof   when   $\nu$ is irrational, we follow
the outline sketched in Section~\ref{sec:irrational}.

 \medskip

As a corollary, we obtain an estimate of the length of {\em any} closed
geodesic along the line of minima $\L_{\mu,\nu}$.
This should be compared with the almost identical estimate on 
p.190 in~\cite{Masur}.
\begin{cor} 
\label{cor:estimateNP}
Let $\mu=a_1\a_1+\dots +a_N\a_N$ and $\nu$ be two measured laminations
which fill up $S$ (where $ \a_1,\dots ,\a_N$ is not necessarily a pants
decomposition). Let  $m_s$ be the minimum of the function $F_s$ and let
$\g$ be  any simple closed curve. Then 
$$
l_{\g}(m_s)=2\sum_{j=1}^N i(\a_j,\g) \log \frac{1}{l_{\a_j}(m_s)}
+o(\log\frac{1}{s}).
$$
\end{cor}
 \begin{proof}
By Theorem~\ref{thm:linetobarycenter}, $m_s\to [\sum\a_j]$. This means
that 
$$
\lim_{s\to 0}\frac{l_{\g}(m_s)}{c_s}
 =i(\g,\textstyle{\sum \a_j}).
$$
 If $\g$ is any closed geodesic, then  $l_{\g}(m_s)=2\sum i(\g,\a_j) 
\log(1/l_{\a_j}(m_s) ) +f(s)$, 
where $f(s)>0$. 
 On the other hand, it is shown in the proof of
Proposition~\ref{prop:estimatedual} that 
$$
\lim \frac{2 \log(1/l_{\a_j}(m_s))}{c_s}=a'_j,
$$
but we know that $a'_j=1$. Therefore  $\lim (f(s)/c_s)=0$ 
   and the result follows.
\end{proof}

\Addresses\recd

\end{document}

%% file: agtout.tex
\def\ifplaintex{\expandafter\ifx\csname documentclass\endcsname\relax}

\def\gtp{{\mathsurround=0pt\it $\cal G\mskip-2mu$eometry \&\ 
$\cal T\!\!$opology $\cal P\!$ublications}}  

\def\recd{{\small Received:\qua\receiveddate\ifx\reviseddate\relax
\else\qquad Revised:\qua\reviseddate\fi\par}} 

\def\lognumber#1{\def\thelognumber{#1}}
\def\volumenumber#1{\def\thevolumenumber{#1}}
\def\volumeyear#1{\def\thevolumeyear{#1}}
\def\papernumber#1{\def\thepapernumber{#1}}
\def\pagenumbers#1#2{\def\startpage{#1}\def\finishpage{#2}}
\def\published#1{\def\publishdate{#1}}

\def\received#1{\def\receiveddate{#1}}

\def\accepted#1{\def\accepteddate{#1}}
\def\asciititle#1{\def\theasciititle{#1}}
\def\covertitle#1{\def\thecovertitle{#1}}
\def\asciiauthors#1{\def\theasciiauthors{#1}}
\def\asciiaddress#1{\def\theasciiaddress{#1}}
\def\asciiemail#1{\def\theasciiemail{#1}}
\def\coverauthors#1{\def\thecoverauthors{#1}}
\long\def\asciiabstract#1{\long\def\theasciiabstract{#1}}
\def\asciikeywords#1{\def\theasciikeywords{#1}}

\let\\\par\let\thelognumber\relax\let\thevolumenumber\relax
\let\thepapernumber\relax\let\thevolumeyear\relax\let\startpage\relax
\let\finishpage\relax\let\publishdate\relax\let\receiveddate\relax
\let\reviseddate\relax\let\accepteddate\relax\let\theasciititle\relax
\let\thecovertitle\relax\let\theasciiauthors\relax\let\theasciiaddress\relax
\let\theasciiabstract\relax\let\theasciikeywords\relax

\let\thecoverauthors\relax\let\theasciiemail\relax

\ifplaintex
\font\logobig=cmssbx10 scaled 3836
\font\logomed=cmssbx10 scaled 2557
\else
\font\logobig=cmssbx10 scaled 4200
\font\logomed=cmssbx10 scaled 2800
\fi

\long\def\makeagttitle{   
\count0=\startpage
\agt\hfill      
\hbox to 45truept{\vbox to 0pt{\vglue -13truept{\logomed A\kern -.37em{\logobig 
T}\kern -.38em G}\vss}\hss}
\break
{\small Volume \thevolumenumber\ (\thevolumeyear)
\startpage--\finishpage\nl
Published: \publishdate}

\vglue .25truein

{\parskip=0pt\leftskip 0pt plus
1fil\def\\{\par\smallskip}{\Large\bf\thetitle}\par\medskip} \vglue
0.05truein

{\parskip=0pt\leftskip 0pt plus 1fil\def\\{\par}{\sc\theauthors}
\par\medskip}%
 
\vglue 0.03truein 

{\small\leftskip 25truept\rightskip 25truept{\bf Abstract}\stdspace\theabstract

{\bf AMS Classification}\stdspace\theprimaryclass
\ifx\thesecondaryclass\relax\else; \thesecondaryclass\fi\par
{\bf Keywords}\stdspace \thekeywords\par}\vglue 7truept

}   

\ifplaintex
\font\phead=cmsl9 scaled 950
\font\pnum=cmbx10 scaled 913
\font\pfoot=cmsl9 scaled 950
\headline{\vbox to 0pt{\vskip -4.5mm\line{\small\phead\ifnum
\count0=\startpage ISSN 1472-2739 (on-line) 1472-2747 (printed)
\hfill {\pnum\folio}\else\ifodd\count0\def\\{ }%
\ifx\theshorttitle\relax\thetitle\else\theshorttitle\fi\hfill{\pnum\folio}
\else\def\\{ and }{\pnum\folio}\hfill\ifx\theshortauthors\relax\theauthors
\else\theshortauthors\fi\fi\fi}\vss}}
\footline{\vbox to 0pt{\vglue 0mm\line{\small\pfoot\ifnum\count0=\startpage
\copyright\ \gtp\hfill\else
\agt, Volume \thevolumenumber\ (\thevolumeyear)\hfill\fi}\vss}}
\else
\font=cmsl9 scaled 1050
\font\lnum=cmbx10 
\font=cmsl9 scaled 1050
\makeatletter
\def\@oddhead{{\small\lhead\ifnum\count0=\startpage ISSN 1472-2739 
(on-line) 1472-2747 (printed)\hfill {\lnum\number\count0}\else\ifodd\count0
\def\\{ }\ifx\theshorttitle\relax \thetitle \else\theshorttitle\fi\hfill
{\lnum\number\count0}\else\def\\{ and }{\lnum\number\count0}
\hfill\ifx\theshortauthors\relax 
\theauthors\else\theshortauthors\fi\fi\fi}}\def\@evenhead{\@oddhead}
\def\@oddfoot{\small\lfoot\ifnum\count0=\startpage\copyright\ \gtp\hfill\else
\agt, Volume \thevolumenumber\ (\thevolumeyear)\hfill\fi}
\def\@evenfoot{\@oddfoot}
\makeatother
\fi
\let\maketitlepage\makeagttitle
\let\makeshorttitle\maketitlepage
\let\maketitle\maketitlepage

\newwrite\gtoutfile
\long\gdef\makeheadfile{  
{\def\\{, }\def\s{ }
\immediate\openout\gtoutfile head.xxx
\immediate\write\gtoutfile{To: math@arxiv.org}
\immediate\write\gtoutfile{Subject: put OR rep NNNNN:ppppp}
\immediate\write\gtoutfile{--text follows this line--}
\immediate\write\gtoutfile{Proxy-for: \ifx\theasciiauthors\relax
\theauthors\else\theasciiauthors\fi\s<\ifx\theasciiemail\relax\theemail\else\theasciiemail\fi>}
\immediate\write\gtoutfile{\noexpand\\}
\immediate\write\gtoutfile{Authors: \ifx\theasciiauthors\relax
\theauthors\else\theasciiauthors\fi}
{\def\\{ }\immediate\write\gtoutfile{Title: \ifx\theasciititle\relax
\thetitle\else\theasciititle\fi}}
\immediate\write\gtoutfile{Subj-class: GT or SG, GR etc}
\immediate\write\gtoutfile{MSC-class: \theprimaryclass\ifx\thesecondaryclass\relax\else, \thesecondaryclass\fi}
\immediate\write\gtoutfile{Journal-ref: Algebr. Geom. Topol. \thevolumenumber\s
(\thevolumeyear) \startpage-\finishpage}
\immediate\write\gtoutfile{Comments: Published by Algebraic and
Geometric Topology at}
\immediate\write\gtoutfile{\s\s\s  http://www.maths.warwick.ac.uk/agt/AGTVol\thevolumenumber/agt-\thevolumenumber-\thepapernumber.abs.html}
\immediate\write\gtoutfile{\noexpand\\}
\immediate\write\gtoutfile{}
\ifx\theasciiabstract\relax
\immediate\write\gtoutfile{\theabstract}\else
\immediate\write\gtoutfile{\theasciiabstract}\fi
\immediate\write\gtoutfile{}
\immediate\write\gtoutfile{\noexpand\\}
\immediate\write\gtoutfile{}
\immediate\closeout\gtoutfile}}  

\def\maketitlepage{\makeagttitle\makeheadfile}
\let\makeshorttitle\maketitlepage
\let\maketitle\maketitlepage